\documentclass[11pt,times,letter]{article}
\usepackage{setspace}

\usepackage{amsmath, amsfonts,amsthm,amssymb, bm}
\usepackage{dsfont}
\usepackage[english]{babel}
\usepackage{float}

\usepackage[normalem]{ulem}

\usepackage{subfigure}
\usepackage{latexsym,amssymb,amsfonts,graphicx}
\usepackage{verbatim}
\usepackage{mathrsfs}
\usepackage{color}
\usepackage{epsfig}
\usepackage{epstopdf}
\usepackage[title]{appendix}
\usepackage{url}
\usepackage{natbib}
\usepackage{tikz}
\usepackage[colorlinks,linkcolor=blue,citecolor=blue,anchorcolor=blue]{hyperref}
\usepackage{mathtools}                              
\mathtoolsset{showonlyrefs=true}          
\newcommand{\idc}{\mathbf{1}}

\newcommand{\Ex}{ \mathbb{E} }

\def\esssup_#1{\underset{#1}{\mathrm{ess\,sup\, }}}
\def\essinf_#1{\underset{#1}{\mathrm{ess\,inf\, }}}
\def\argmax_#1{\underset{#1}{\mathrm{arg\,max\, }}}
\def\argmin_#1{\underset{#1}{\mathrm{arg\,min\, }}}

\newcommand{\R}{\mathds{R}}

\newtheorem{theorem}{Theorem}[section]

\numberwithin{equation}{section}

\newtheorem{proposition}[theorem]{Proposition}
\newtheorem{remark}[theorem]{Remark}
\newtheorem{lemma}[theorem]{Lemma}
\newtheorem{corollary}[theorem]{Corollary}
\newtheorem{ansatz}[theorem]{Ansatz}

\allowdisplaybreaks

\definecolor{Red}{rgb}{1.00, 0.00, 0.00}

\definecolor{DRed}{rgb}{0.5, 0.00, 0.00}

\definecolor{Blue}{rgb}{0.00, 0.00, 1.00}

\definecolor{Green}{rgb}{0.0, 0.4, 0.0}

\title{Optimal Annuitization, Consumption, and Investment with Costly Capital Injections}
\author{Yijie Huang%
\thanks{Department of Applied Mathematics, The Hong Kong Polytechnic University,
        Kowloon, Hong Kong, China. Email: yijie.huang@polyu.edu.hk.}
\and Xiaoqing Liang%
\thanks{Corresponding author. School of Sciences, Hebei University of Technology, Tianjin, 300401, P. R. China. Email: liangxiaoqing115@hotmail.com.}
\and Jie Xiong%
\thanks{Department of Mathematics and SUSTech International center for Mathematics, Southern University
			of Science and Technology, Shenzhen, 518055, P. R. China. Email: xiongj@sustech.edu.cn.}
\and Xiang Yu%
\thanks{Department of Applied Mathematics, The Hong Kong Polytechnic University,
       Kowloon, Hong Kong, China. Email: xiang.yu@polyu.edu.hk.}
	}

\date{\vspace{-0.5in}}

\begin{document}
\maketitle

\begin{abstract}
We study a joint optimal decision of portfolio, consumption, and annuitization for a pension-plan manager who balances retirement consumption benefits against the cost of maintaining the solvency of pension account. The pension sponsor provides costly capital injection to prevent negative wealth. Annuity purchases are irreversible, and injected funds remain in the liquid account and earn interest.  Using the Legendre transform, we derive a semi-explicit solution to the associated free-boundary problem and verify the optimality of the resulting strategies.  We identify an endogenous annuity-income threshold separating injection-active and injection-inactive regions and characterize six policy regimes arising from the interplay between capital injections and consumption constraints.  Compared with the no-injection benchmark, capital injections delay annuitization, alter the structure of the annuitization boundary, and permit positive holdings in risky assets even at zero liquid wealth. Numerical results further reveal that cheaper annuities may increase the expected capital contributions required from the pension sponsor.
\end{abstract}


\section{Introduction}\label{sec:intro}
After the seminal work \citet{Merton69,Merton71}, continuous-time utility maximization for portfolio management and consumption planning has become a core research topic in quantitative finance. Taking into account the retirement planning, these decisions become more subtle because of longevity risk and the availability of life annuities, which convert liquid wealth into a guaranteed income stream for the remaining lifetime.  As shown by \citet{Yaari65}, actuarially fair annuities can be preferable to bonds in the absence of a bequest motive, while \citet{DavidoffBrownDiamond05} establish the welfare benefits of annuitization under more general conditions. These findings highlight the importance of jointly optimizing consumption, investment, and annuitization decisions in retirement planning.

A stream of studies investigated dynamic annuitization decisions when retirees can invest in financial markets and choose when to convert liquid wealth into lifetime income. The irreversibility of annuity purchases creates a trade-off between guaranteed lifetime income and the flexibility of retaining liquid wealth, particularly under non-negative consumption and wealth constraints. \citet{MilevskyMooreYoung06} studied optimal annuitization and investment strategies aimed at minimizing lifetime financial ruin, while \citet{MilevskyYoung07} examined optimal consumption, investment, and annuitization under different annuity-market arrangements. \citet{WangYoung12} extended this framework to commutable annuities with surrender charges. Under CARA utility, \citet{LiangYoung18} imposed non-negative consumption and wealth and solve a two-dimensional free-boundary problem by Legendre duality; see also \cite{DLPY22} for the exponential utility on consumption with reference to past spending maximum.  CARA measures risk tolerance in dollars but lacks the homogeneity of power utility, making the constraints studied by \citet{KaratzasEtAl86} and \citet{CoxHuang89} especially relevant.  More recent extensions consider deferred annuities and refund provisions \citep{HabibEtAl20}, housing and means-tested pensions \citep{AndreassonShevchenko24}, stochastic interest rates \citep{DillschneiderEtAl24}, and fixed transaction costs \cite{GLSX2026}.

Downside protection is another important consideration in post-retirement pension management and annuitization, as adverse investment outcomes and the depletion of liquid wealth may threaten the sustainability of retirement consumption. In income-drawdown arrangements, retirees may remain invested during the decumulation phase rather than fully annuitizing their pension wealth, making protection against financial shortfalls particularly relevant. \citet{GerrardHabermanVigna04} studied optimal post-retirement investment, while \citet{GerrardHojgaardVigna12} jointly examined consumption, investment, and annuitization timing. More explicitly incorporating downside protection, \citet{DiGiacintoFedericoGozziVigna14} investigated income drawdown with a minimum guarantee, while \citet{Dadashi20} considered optimal investment and consumption under a minimum guarantee for future annuity income. These studies highlight the importance of balancing investment flexibility and annuitization decisions with protection against the risk of insufficient retirement wealth. However, they generally impose non-negative wealth as a hard state constraint or provide minimum retirement benefits, without explicitly modeling the capital contributions required to maintain solvency. In practice, pension sponsors may be required to make additional contributions when funding levels deteriorate, as illustrated by funding requirements for U.S. defined-benefit pension plans \citep{Rauh06}.

The provision of such external financial support is closely related to the literature on costly capital injections and reflected surplus processes. In dividend and insurance models, capital injections are used to prevent insolvency, often leading to a lower reflecting boundary at zero \citep{LokkaZervos08,EisenbergSchmidli09,EisenbergSchmidli11,FerrariSchuhmann19}. Related portfolio optimization with benchmark-tracking \citep{BoLiaoYu21,BoHuangYu24,BoHuangYu25SICON,BoHuangYu25,BHYY26} also utilized the capital support to maintain prescribed wealth floors, with the minimal regulator representing cumulative shortfalls. These studies provided a mathematical framework for modeling the cost of maintaining solvency through external funding. However, unlike dividend payments, which distribute surplus without changing the underlying income structure, annuity purchases convert liquid wealth into a flow of lifetime income, thereby affecting both current wealth and future cash flows. Incorporating irreversible annuitization alongside costly capital injections thus introduces an additional dimension to the stochastic control problem.

Building on these two strands of literature, we consider a post-retirement pension account with a sponsor-backed solvency guarantee. The pension-plan manager jointly determines consumption, portfolio allocation, and irreversible annuity purchases, aiming to maximize the retiree's expected discounted consumption utility under constant absolute risk aversion (CARA) while accounting for the cost of external capital support. Whenever adverse investment outcomes or consumption withdrawals would otherwise drive liquid wealth below zero, the pension sponsor provides the minimum amount of capital required to maintain solvency. The injected funds remain in the pension account and earn interest with remaining wealth. Importantly, capital injections are available solely to fulfill the solvency guarantee and cannot be used as a discretionary borrowing facility. Unlike annuity purchases, which constitute a singular control chosen by the manager, capital injections are determined endogenously as the minimal reflecting process required to maintain non-negative liquid wealth. The resulting problem therefore involves a lower reflecting boundary associated with solvency protection and an upper intervention boundary governing irreversible annuity purchases.

To solve this problem, we first derive the associated Hamilton-Jacobi-Bellman variational inequality (HJB-VI), which incorporates a gradient constraint arising from irreversible annuity purchases, and a lower-boundary marginal-value constraint coming from costly capital injections. Unlike conventional annuitization problems with a hard non-negativity constraint on wealth, the lower-boundary constraint may be binding or non-binding depending on the existing annuity income, leading to an endogenous transition between injection-active and injection-inactive regions. The interaction between capital injections at the lower wealth boundary and annuity purchases at an upper intervention boundary gives rise to a two-dimensional free-boundary problem. Using the Legendre transform, we derive a semi-explicit solution and establish the characterization of the optimal strategies via the verification theorem. In particular, we identify an endogenous annuity-income threshold that separates the injection-active and injection-inactive regions. Below this threshold, the pension sponsor may need to provide capital support to maintain solvency, whereas above it, the optimal investment and consumption policies preserve non-negative wealth without further injections. Together with the non-negativity constraint on consumption, this distinction gives rise to six policy regimes, characterized by active or inactive capital injections and zero, mixed, or positive consumption. We further characterize the corresponding optimal investment, consumption, and annuitization strategies, revealing how costly solvency support affects the allocation between liquid wealth and guaranteed lifetime income.

We further conduct numerical experiments to illustrate the financial implications of the optimal strategies and the interplay between annuitization and solvency protection. Compared with the benchmark model in \citet{LiangYoung18}, the availability of capital injection raises the optimal annuitization boundary, encouraging the manager to retain liquid wealth for longer before purchasing additional annuities. Moreover, the annuitization boundary decreases with existing annuity income in the injection-active region, in contrast to the increasing boundary in the no-injection model. Capital injections also allow the manager to maintain positive risky investment positions even at zero liquid wealth, highlighting the role of external solvency support in facilitating risk taking. Our sensitivity analysis reveals two new observations. First, although cheaper annuities encourage earlier annuitization, they may increase the expected discounted capital contributions required from the pension sponsor. Second, the optimal annuitization boundary may respond non-monotonically to subjective mortality when capital injections are available. These findings provide new insights in the joint decision of annuitization, portfolio risk taking, and the funding costs of solvency guarantees in pension-plan management.

The remainder of this paper is organized as follows. Section~\ref{sec:model} introduces the financial market and formulates the pension management problem with costly capital injections. Section~\ref{sec:dual} derives a semi-explicit solution to the associated HJB variational inequality using the Legendre transform, while Section~\ref{sec:verification} establishes the optimality of the resulting strategies through a verification theorem. Section~\ref{sec:numerical} presents numerical examples illustrating the optimal policies, their comparison with the no-injection benchmark, and the economic implications of annuitization and solvency protection. Appendix \ref{sec:special-cases} presents the result for the special case when the interest rate $r=0$ and Appendix \ref{sec:proof} collects all the proofs of the main results in the paper.

\section{Problem Formulation}\label{sec:model}
We consider a retiree participating in a pension plan that provides solvency protection for her retirement account. The plan commits to maintaining the account at a non-negative liquid wealth level by supplying additional capital whenever the account would otherwise become insolvent. Meanwhile, the pension-plan manager determines the payout rate available for retirement consumption, the portfolio allocation, and the annuitization policy. 

\subsection{Financial and Annuity Markets}
Let $(\Omega,\mathcal F,\mathbb F,\mathbb P)$ be a filtered probability space satisfying the usual conditions, where $\mathbb F=(\mathcal F_t)_{t\geq0}$ supports a one-dimensional Brownian motion $W$.  The pension-plan manager can trade continuously in one riskless and one risky asset, whose prices follow the dynamics
\begin{align*}
    \frac{dB_t}{B_t}=r\,dt,\quad B_0=1,\qquad
    \frac{dS_t}{S_t}=\mu\,dt+\sigma\,dW_t,\quad S_0>0,
\end{align*}
where the risk-free rate $r> 0$, the return rate $\mu>r$ and the volatility $\sigma>0$ are constants.  

Let $\tau_d$ denote the retiree's remaining lifetime.  We assume that $\tau_d$ is independent of $W$ and exponentially distributed with subjective mortality intensity $\lambda^s>0$; hence $\mathbb P(\tau_d>t)=e^{-\lambda^s t}$, for $t\geq0$. The intensity $\lambda^s$ represents the retiree's survival belief and therefore determines the effective horizon in her objective.

At any time, the pension-plan manager may irreversibly purchase a single-premium immediate life annuity.  One unit of annuity income pays continuously at rate one until death and has no surrender value or bequest value.  Insurers price this income using the risk-free rate and a constant objective mortality intensity $\lambda^o>0$.  Its single premium is
\begin{align}\label{eq:theta}
    \theta = \int_0^\infty e^{-(r+\lambda^o)t}\,dt
    =\frac{1}{r+\lambda^o}.
\end{align}
Note that $\lambda^o$ may differ from $\lambda^s$, allowing the insurer's pricing basis to diverge from the retiree's survival belief and also capturing a reduced-form loading for the annuity \citep{LiangYoung18}.

\subsection{Pension Plan with a Solvency Guarantee}

At time $t$, let $\pi_t\in\mathbb R$ denote the amount invested in the risky asset and let $c_t\geq0$ denote the payout rate available for retirement consumption chosen by the pension-plan manager. The process $A=(A_t)_{t\geq 0}$ represents the retiree's total annual annuity income, including both pre-existing annuity income and additional annuities purchased by the plan up to time $t$. Let $X_t$ denote the liquid wealth in the retiree's pension account. Capital injected by the pension plan enters the funded account and, once injected, is invested together with the remaining liquid wealth. Under the strategy $(\pi,c,A)$, the liquid wealth process evolves according to the following dynamics:
\begin{align}\label{state-X}
dX_t=&\left(rX_{t-}+(\mu-r)\pi_{t-}-c_{t-}+A_{t-}\right)dt+\sigma\pi_{t-}dW_t-\theta dA_t+dL_t^X, \quad X_{0-}=x\geq0, 
\end{align}
 where $L^X$ denotes cumulative capital injected by the pension plan. 
 
 The pension plan guarantees that the pension account remains solvent at all times. Accordingly, $L^X$ is the minimal non-decreasing process required to enforce the solvency constraint $X_t\geq 0$ with $L_0^X=0$ and
$\int_0^\infty {\bf 1}_{\{X_t>0\}}dL_t^X=0$. Thus, capital is supplied only when the funded account reaches the solvency boundary and only in the amount needed to prevent negative wealth. It is important to emphasize that $L^X$ is not an additional discretionary source of borrowing available to the manager; rather, it represents the minimum capital contribution required for the pension plan to meet its solvency guarantee.

\subsection{Manager's Optimization Problem}\label{sec:optimization}
The pension-plan manager jointly determines consumption, investment, and annuitization strategies to balance the retiree’s utility from consumption against the expected future cost of capital injections required to maintain the plan’s solvency guarantee.  The admissible regular and annuitization controls in $\mathbb U$ satisfy the following conditions on every finite horizon $T$:
\begin{itemize}
    \item The control processes $\pi=(\pi_t)_{t\geq 0}$ and $c=(c_t)_{t\geq 0}$ taking values on  $\R\times\R_+$ ($\R_+:=[0,\infty)$) are progressively measurable, and satisfy
$\Ex[\int_0^T(\pi_t^2+c_t)dt]<\infty$.
    \item The annuitization strategy $A=(A_t)_{t\geq 0}$ is an adapted, non-negative, non-decreasing process whose paths are right-continuous with left limits.
\end{itemize}

We assume that the retiree's preference is modelled by the exponential utility $U(x) = -\frac{1}{\gamma} e^{-\gamma x}$, where $\gamma>0$ is absolute risk aversion. The non-negative consumption constraint becomes active under the choice of exponential utility. The pension-plan manager's objective is to optimize the  trade-off between the retiree's utility from consumption and the cost of capital injections borne by the plan. The corresponding value function is defined by
\begin{align}\label{eq:u0}
V(x,a):=\sup_{(\pi,c,A)\in\mathbb U}\Ex\left[\int_0^{\tau_d} e^{-\rho t}U(c_t)dt-\beta\int_0^{\tau_d} e^{-\rho t}dL_t^X \Big|X_{0-}=x,A_{0-}=a\right],
\end{align}
where $\rho>0$ is the discount rate and $\beta>0$ is the marginal utility-equivalent cost to the pension plan  of injecting one unit of capital.  By independence of lifetime $\tau_d$ and the Brownian motion $W$, we can rewrite the value function given by \eqref{eq:u0} as follows,
\begin{align}\label{eq:u1}
V(x,a)=\sup_{(\pi,c,A)\in\mathbb U}\Ex\left[\int_0^\infty e^{-\delta t}U(c_t)dt-\beta\int_0^\infty e^{-\delta t}dL_t^X \Big|X_{0-}=x,A_{0-}=a\right],
\end{align}
where $\delta:=\rho+\lambda^s>0$ is the effective discount rate that incorporates the subjective mortality intensity $\lambda^s$.  The first term in \eqref{eq:u1} measures the retiree's discounted lifetime utility from consumption, whereas the second term captures the discounted cost of capital injections associated with the pension plan's solvency guarantee.

\subsection{HJB Variational Inequality }
Applying dynamic programming to \eqref{eq:u1} yields the HJB variational inequality (HJB-VI), including the lower-boundary gradient constraint, for $(x,a)\in\R_+^2$:
\begin{align}\label{eq:HJB-v}
\begin{cases}
\displaystyle\max\bigg\{\sup_{\pi \in\R}\left((\mu-r)\pi v_x+\frac{1}{2}\sigma^2 \pi^2 v_{xx}\right)+\sup_{c\geq 0}\left(U(c)-cv_x\right)+(rx+a)v_x-\delta v,v_a-\theta v_x \bigg\}=0,\\[1em]
\displaystyle v_x(0,a)\leq\beta,\quad \forall a\in\R_+.
\end{cases}
\end{align}
The second entry in \eqref{eq:HJB-v}  is the directional gain from an infinitesimal annuity purchase: one unit of new income raises $a$ by one and reduces liquid wealth by $\theta$.  Hence $v_a-\theta v_x<0$ in the no-purchase region and equals zero when annuitization is optimal. The lower-boundary inequality has a different origin.  One unit of reflection increases wealth by one at marginal cost $\beta$, so optimality requires $v_x(0,a)\leq\beta$.  If capital injection is active, complementarity gives the Neumann condition $v_x(0,a)=\beta$.  If the optimal regular controls make the boundary viable without reflection, then $v_x(0,a)<\beta$ and the hard state-constraint conditions $\pi^*(0,a)=0$ and $c^*(0,a)\leq a$ apply.  In either case,
\begin{align}\label{eq:injection-complementarity}
    \big(\beta-v_x(0,A_t^*)\big)dL_t^{X^*}=0.
\end{align}
This condition is the precise switch between the injection-active and injection-inactive solutions developed below.

Motivated by \citet{LiangYoung18}, we make the following ansatz for the
optimal strategy.

\begin{ansatz}[Two-sided barrier policy]\label{ans:two-sided-barrier}
For each current annuity income $a$, there is an annuitization boundary
$b(a)>0$.  No annuity is bought while $0\leq x<b(a)$.  When the state reaches
$b(a)$, annuity income is purchased just fast enough to keep the state at 
that boundary; if $x>b(a)$ initially, the purchase $\Delta a$ is determined by
\begin{align}\label{eq:projection}
    x-\theta\Delta a=b(a+\Delta a).
\end{align}
At the lower boundary, capital injection is used only when the marginal-value
cap binds.  Thus there is a switching income $a_\beta$ such that
$v_x(0,a)=\beta$ for $a<a_\beta$, whereas for $a>a_\beta$ the hard
state-constraint conditions $\pi^*(0,a)=0$ and $c^*(0,a)\leq a$ make the
boundary viable without injection.  Finally, because $U'(0)=1$, consumption is
zero when $v_x\geq1$ and positive when $v_x<1$.
\end{ansatz}

Under this ansatz, we compute the corresponding expected utility function $v$ in each region subject to capital injection and determine
the two dual boundaries; then recover the candidate value function and verify the optimality of the proposed strategies in Section~\ref{sec:verification}.  

\section{Solving for the Value Function}\label{sec:dual}

In this section, we solve the HJB-VI given in \eqref{eq:HJB-v} for the value function. Under Ansatz~\ref{ans:two-sided-barrier}, define the continuation and annuitization
regions as
\begin{align}\label{eq:regions}
    \mathcal C=\{(x,a):0\leq x<b(a)\},
    \qquad
    \mathcal A=\{(x,a):x\geq b(a)\}.
\end{align}
On the continuation region $\mathcal C$, the candidate function $v$ solves the following free-boundary problem:
\begin{align}\label{eq:free-boundary-prob}
\sup_{\pi \in\R}\left((\mu-r)\pi v_x+\frac{1}{2}\sigma^2 \pi^2 v_{xx}\right)+\sup_{c\geq 0}\left(U(c)-cv_x\right)+(rx+a)v_x-\delta v = 0,
\end{align}
subject to the boundary condition $v_x(0,a)=\beta$ if capital injection is active. Otherwise, if the optimal regular controls make the boundary viable without reflection, then $v_x(0,a)<\beta$, and the following the lower state-constraint conditions hold:
\begin{align}\label{eq:lower_state_constraint}
\pi(0,a)=0 \quad~\mbox{and}~\quad  0 \le c(0,a)\leq a.
\end{align}

Suppose that $v_x>0$ and $v_{xx}<0$, and then the first-order conditions in the no-annuitization region give
\begin{align}\label{eq:feedback-regular}
\pi^*(x,a)&=-\frac{\mu-r}{\sigma^2}\frac{v_x(x,a)}{v_{xx}(x,a)},
&
c^*(x,a)&=\left(-\frac{1}{\gamma}\log v_x(x,a)\right)^+.
\end{align}
Let $m:=\frac{(\mu-r)^2}{2\sigma^2}>0$. Define
\begin{align}\label{eq:utility-conjugate}
\widetilde U(y)&:=\sup_{c\geq0}\{U(c)-cy\}
=\begin{cases}
\dfrac{y\log y-y}{\gamma},&0<y<1,\\[0.8em]
-\dfrac{1}{\gamma},&y\geq1.
\end{cases}
\end{align}
The free boundary problem \eqref{eq:free-boundary-prob} therefore reduces to
\begin{align}\label{eq:primal-continuation}
m\frac{v_x^2}{v_{xx}}+\delta v-(rx+a)v_x-\widetilde U(v_x)=0.
\end{align}
On the purchase-annuity region $\mathcal A$, $v$ satisfies
$- \theta v_x(x, a) + v_{a}(x, a) = 0$. Particularly, on the annuitization boundary $x = b(a)$, we have
$- \theta v_x(b(a), a) + v_{a}(b(a), a) = 0$.

\subsection{Legendre transform}
For $y\in(0,\beta]$, define the dual function  $\widehat v(y,a):=\sup_{x\geq0}\{v(x,a)-xy\}$. At an interior maximizer, we have
\begin{align}\label{eq:dual-relations}
    y=v_x(x,a),\quad x=-\widehat v_y(y,a),\quad
    v(x,a)=\widehat v(y,a)-y\widehat v_y(y,a),
    \quad v_{xx}(x,a)=-\frac{1}{\widehat v_{yy}(y,a)}.
\end{align}
Substituting \eqref{eq:dual-relations} into \eqref{eq:primal-continuation} and collecting terms, we obtain:
\begin{align}\label{eq:dual-ode}
    m y^2\widehat v_{yy}+(\delta-r)y\widehat v_y-\delta\widehat v
    +ay+\widetilde U(y)=0.
\end{align}
This is an Euler equation in $y$, with $a$ acting as a parameter between annuity purchases.  Its characteristic roots are
\begin{align}\label{eq:general-roots-v+} 
\nu_{\pm}
&=\frac{m-\delta+r \pm
\sqrt{(\delta-r-m)^2+4m\delta}}{2m},  
\end{align}
with $\nu_+>1$ and $\nu_-<0$.
The solutions on the two consumption branches are
\begin{align}\label{eq:dual-general-r-low}
\widehat v_+(y,a)
={}&C_1(a)y^{\nu_+}+C_2(a)y^{\nu_-}
+\frac{1}{\gamma r}y\log y
+\left(\frac{a}{r}+\frac{m+\delta-2r}{\gamma r^2}\right)y,
\qquad 0<y\leq1,
\end{align}
and
\begin{align}\label{eq:dual-general-r-high}
\widehat v_0(y,a)
={}&C_3(a)y^{\nu_+}+C_4(a)y^{\nu_-}
+\frac{a}{r}y-\frac{1}{\gamma\delta},
\qquad y\geq1.
\end{align}
The subscripts $+$ and $0$ of $\hat{v}$ denote positive and zero consumption, respectively. Evaluating the two branches and their first derivatives at $y=1$ gives
\begin{align*}
& [C_1(a)-C_3(a)]+[C_2(a)-C_4(a)]
 =-\frac{1}{\gamma\delta}
   -\frac{m+\delta-2r}{\gamma r^2},\\
& \nu_+[C_1(a)-C_3(a)]+\nu_-[C_2(a)-C_4(a)]
 =-\frac{m+\delta-r}{\gamma r^2}.
\end{align*}
Solving the above equations, we obtain
\begin{align}\label{eq:general-r-matching}
C_1(a)-C_3(a)&=\Delta_+,
&
C_2(a)-C_4(a)&=\Delta_-,
\end{align}
in which
\begin{align}
\Delta_+:=
\frac{-\dfrac{m+\delta-r}{\gamma r^2}
+\nu_-\left(\dfrac{1}{\gamma\delta}
+\dfrac{m+\delta-2r}{\gamma r^2}\right)}{\nu_+-\nu_-},\quad
\Delta_-:=
\frac{\dfrac{m+\delta-r}{\gamma r^2}
-\nu_+\left(\dfrac{1}{\gamma\delta}
+\dfrac{m+\delta-2r}{\gamma r^2}\right)}{\nu_+-\nu_-}. \\ \label{eq:general-r-matching-data+1}
\end{align}
Using \eqref{eq:general-roots-v+} , we rewrite \eqref{eq:general-r-matching-data+1} as
\begin{align}\label{eq:general-r-matching-data+}
\Delta_+&
=-\frac{1-\nu_-}
 {\gamma r\nu_+(\nu_+-1)(\nu_+-\nu_-)}<0 ,
\quad
\Delta_-
=\frac{\nu_+-1}{\gamma r\nu_-(1-\nu_-)(\nu_+-\nu_-)}<0.
\end{align}
Noting that the differences in \eqref{eq:general-r-matching} do not depend
on $a$, we obtain $C_1'(a)=C_3'(a)$ and $C_2'(a)=C_4'(a)$.

Let the upper and lower endpoints of the continuation interval in the dual variable be
\begin{align}\label{eq:dual-annuity-boundary}
    y_0(a):=v_x(0,a),\qquad y_b(a):=v_x(b(a),a),
    \qquad 0<y_b(a)<y_0(a)\leq\beta.
\end{align}
The boundary conditions become
\begin{align}
\widehat v_y(y_0(a),a)&=0, \label{eq:dual-boundaries}\\
\widehat v_a(y_b(a),a)&=\theta y_b(a), \qquad
\widehat v_{ay}(y_b(a),a)=\theta. \label{eq:dual-boundaries-annuitization-condition1}
\end{align}
When capital injection is inactive, the zero-volatility state constraint $\pi^*(0,a)=0$ gives the additional super-contact condition
\begin{align}\label{eq:dual-state-constraint}
    \widehat v_{yy}(y_0(a),a)=0.
\end{align}
When capital injection is active, the super-contact condition is replaced
by $y_0(a)=\beta$.  Thus, the two lower-boundary cases are different:
in the inactive case, the endpoint is unknown and the curvature vanishes
there; in the active case, the endpoint is prescribed and the curvature is
to be determined.

\subsection{Determining the dual boundaries}\label{sec:positive-rate-construction}
We first use the annuitization conditions to determine the derivatives of the
coefficients $C_j(a), j=1,2,3,4$. The condition at zero wealth will then determine the ratio
of the endpoints $y_0(a)$ and $y_b(a)$.  After that, we solve for the endpoints themselves.

\paragraph{The annuitization conditions.}
Both particular solutions of \eqref{eq:dual-general-r-low} and \eqref{eq:dual-general-r-high} have $a$-derivative $y/r$.  Moreover,
$C_1'=C_3'$ and $C_2'=C_4'$.  Thus the two conditions at $y_b(a)$ give
the same equations on either consumption branch:
\begin{align}
& C_1'(a)y_b(a)^{\nu_+}+C_2'(a)y_b(a)^{\nu_-}
 =-\frac{1-r\theta}{r}y_b(a),\label{eq:coefficient-derivative-1}\\
& \nu_+C_1'(a)y_b(a)^{\nu_+-1}
 +\nu_-C_2'(a)y_b(a)^{\nu_--1}
 =-\frac{1-r\theta}{r}. \label{eq:coefficient-derivative-2}
\end{align}
Eliminating $C_1'(a)$ and $C_2'(a)$ from the two equations above gives
\begin{align}
& C_1'(a)=C_3'(a)
 =-\frac{1-r\theta}{r}\frac{1-\nu_-}{\nu_+-\nu_-}
       y_b(a)^{1-\nu_+}<0,\label{eq:coefficient-derivatives13}\\
& C_2'(a)=C_4'(a)
=-\frac{1-r\theta}{r}\frac{\nu_+-1}{\nu_+-\nu_-}
       y_b(a)^{1-\nu_-}<0.\label{eq:coefficient-derivatives24}
\end{align}
 Substituting these
expressions into $\widehat v_{ay}$ gives
\begin{align}\label{eq:annuity-mixed-derivative}
 \widehat v_{ay}(y,a)
 =\frac1r-\frac{1-r\theta}{r}
 \left[
 \frac{\nu_+(1-\nu_-)}{\nu_+-\nu_-}
       \left(\frac{y}{y_b(a)}\right)^{\nu_+-1}
 +\frac{\nu_-(\nu_+-1)}{\nu_+-\nu_-}
       \left(\frac{y}{y_b(a)}\right)^{\nu_--1}
 \right].
\end{align}

\paragraph{The condition at zero wealth.}
Differentiating $\widehat v_y(y_0(a),a)=0$ with respect to $a$ along the boundary gives
\[
 \widehat v_{ay}(y_0(a),a)
 +\widehat v_{yy}(y_0(a),a)y_0'(a)=0.
\]
If injection is inactive, the second term vanishes by super-contact.
If injection is active, it vanishes because $y_0(a)=\beta$ is constant.
Thus, in both cases, 
\begin{align}\label{eq:deriv_hv_ay}
\widehat v_{ay}(y_0(a),a)=0.
\end{align}
Substitution into \eqref{eq:annuity-mixed-derivative} gives an equation for the ratio
$z_{0b}:=y_0(a)/y_b(a)$:
\begin{align}\label{eq:z0b-general}
 (1-r\theta)\left[
 \frac{\nu_+(1-\nu_-)}{\nu_+-\nu_-}z_{0b}^{\nu_+-1}
 +\frac{\nu_-(\nu_+-1)}{\nu_+-\nu_-}z_{0b}^{\nu_--1}
 \right]=1.
\end{align}
The equation is independent of $a$.  Define
\begin{align}\label{eq:kappa-general}
 \kappa:=\frac{1-\nu_-(1-r\theta)z_{0b}^{\nu_--1}}{1-\nu_-}.
\end{align}

\begin{lemma}\label{lem:boundary-ratio}
Equation \eqref{eq:z0b-general} has exactly one solution $z_{0b}>1$.
Consequently, $y_0(a)/y_b(a)$ is independent of $a$, for both lower-boundary
alternatives.  The constant $\kappa$ satisfies $\frac1{1-\nu_-}<\kappa<1$.
\end{lemma}

Equation \eqref{eq:z0b-general} and the definition of $\kappa$ also give
\begin{align}\label{eq:kappa-root-identity1}
 (1-r\theta)\nu_+z_{0b}^{\nu_+-1}
 =1+\kappa(\nu_+-1),\quad
 (1-r\theta)\nu_-z_{0b}^{\nu_--1}
 =1+\kappa(\nu_--1). 
\end{align}
Differentiating \eqref{eq:annuity-mixed-derivative} once more in $y$,
setting $y=y_0(a)$, and using \eqref{eq:z0b-general}, \eqref{eq:kappa-root-identity1}, and
$m(\nu_+-1)(1-\nu_-)=r$ yields
\begin{align}
 \widehat v_{ayy}(y_0(a),a) =-\frac{(\nu_+-1)
 [1-\nu_-(1-r\theta)z_{0b}^{\nu_--1}]}{r\,y_0(a)}
 =-\frac{\kappa}{m y_0(a)}. \label{eq:endpoint-mixed-curvature}
\end{align}

For the annuity-purchase inequality, define
$g(y,a):=\widehat v_a(y,a)-\theta y$.  Differentiation of the dual equation \eqref{eq:dual-ode} with respect to $a$ gives
\begin{align}\label{eq:general-r-g}
 my^2g_{yy}+(\delta-r)yg_y-\delta g+(1-r\theta)y=0,
 \qquad g(y_b(a),a)=g_y(y_b(a),a)=0.
\end{align}
Here the coefficients of the homogeneous terms are already determined
by \eqref{eq:coefficient-derivatives13}--\eqref{eq:coefficient-derivatives24}.  Hence
\begin{align}\label{eq:general-r-g-explicit}
\begin{split}
 g(y,a)=\frac{1-r\theta}{r}y_b(a)
 \left[
 \frac{y}{y_b(a)}
 -\frac{1-\nu_-}{\nu_+-\nu_-}
       \left(\frac{y}{y_b(a)}\right)^{\nu_+}
 -\frac{\nu_+-1}{\nu_+-\nu_-}
       \left(\frac{y}{y_b(a)}\right)^{\nu_-}
 \right].
\end{split}
\end{align}
The negative power coefficients and
$\nu_\pm(\nu_\pm-1)>0$ imply $g_{yy}<0$.
Since $g=g_y=0$ at $y_b(a)$, we obtain $g<0$ for $y>y_b(a)$.

\subsection{Solving the no-injection problem}
We use the superscript $\mathrm{N}$ for the problem in which no capital
injection is allowed.  Let $y_0^{\rm N}(a)$ denote its unknown
zero-wealth endpoint and set $y_b^{\rm N}(a)=y_0^{\rm N}(a)/z_{0b}$.
The zero-wealth and super-contact conditions, respectively, are
\begin{align}\label{eq:positive-rate-contact}
 \widehat v_y^{\rm N}(y_0^{\rm N}(a),a)=0,\qquad
 \widehat v_{yy}^{\rm N}(y_0^{\rm N}(a),a)=0.
\end{align}
We first solve the positive-consumption case, then continue the solution
through the mixed case to the zero-consumption case. The annuity income thresholds are determined as part of the analysis.

\subsubsection{Positive consumption: $a \ge \kappa/\gamma$}
Suppose first that $e^{- a \gamma} \le y_0^{\rm N}(a)\leq1$, so the consumption rate in whole continuation
interval is positive.  Differentiating
the dual equation \eqref{eq:dual-ode} in $y$ gives
\[
 my^2\widehat v_{yyy}^{\rm N}
 +(2m+\delta-r)y\widehat v_{yy}^{\rm N}
 -r\widehat v_y^{\rm N}
 +a+\widetilde U'(y)=0.
\]
At the zero-wealth endpoint, both the first and second derivatives
vanish.  It follows that
\begin{align}\label{eq:endpoint-third-derivative}
 \widehat v_{yyy}^{\rm N}(y_0^{\rm N}(a),a)
 =-\frac{a+\widetilde U'(y_0^{\rm N}(a))}
         {my_0^{\rm N}(a)^2}.
\end{align}
Next differentiate the super-contact condition
$\widehat v_{yy}^{\rm N}(y_0^{\rm N}(a),a)=0$ with respect to $a$ along the endpoint
$y=y_0^{\rm N}(a)$.  Using 
\eqref{eq:endpoint-mixed-curvature} and \eqref{eq:endpoint-third-derivative} gives
\[
 0=-\frac{\kappa}{my_0^{\rm N}(a)}
 -\frac{a+\widetilde U'(y_0^{\rm N}(a))}
        {my_0^{\rm N}(a)^2}(y_0^{\rm N})'(a).
\]
Thus, $y_0^{\rm N}(a)$ satisfies the following equation
\begin{align}\label{eq:uncapped-endpoint-equation}
 \left(a+\widetilde U'(y_0^{\rm N}(a))\right)
 \frac{(y_0^{\rm N})'(a)}{y_0^{\rm N}(a)}=-\kappa.
\end{align}
Since in the present case $\widetilde U'(y)=\log y/\gamma$,  multiplication
by the integrating factor $y_0^{\rm N}(a)^{1/\kappa}$ shows that
\[
 \frac{d}{da}\left\{
 y_0^{\rm N}(a)^{1/\kappa}
 \left[a+\frac{\log y_0^{\rm N}(a)-\kappa}{\gamma}\right]
 \right\}=0.
\]
Let $a_2$ denote the annuity income at which
$y_0^{\rm N}(a_2)=1$.  Evaluating the constant of integration there
gives
\begin{align}\label{eq:positive-endpoint-integrated}
 a=\left(a_2-\frac{\kappa}{\gamma}\right)
       y_0^{\rm N}(a)^{-1/\kappa}
   +\frac{\kappa-\log y_0^{\rm N}(a)}{\gamma}.
\end{align}
However \eqref{eq:positive-endpoint-integrated} alone does not determine $a_2$. We determine $a_2$ below using the bounds on
the value function.

To select the solution compatible with the objective function, note that a
no-injection value function must satisfy
\[
 -\frac{e^{-\gamma a}}{\gamma\delta}\leq v^{\rm N}(0,a)\leq0,
\]
where the lower bound is derived by consuming the existing annuity income
without investment or further purchases.
At zero wealth, evaluate the dual equation \eqref{eq:dual-ode} at
$y=y_0^{\rm N}(a)$. Using the positive-consumption branch of
$\widetilde U$ in \eqref{eq:utility-conjugate}, the two contact conditions in
\eqref{eq:positive-rate-contact}, and the dual relationship
$\widehat v^{\rm N}(y_0^{\rm N}(a),a)=v^{\rm N}(0,a)$, we obtain
\begin{align}\label{eq:inactive-lower-value}
 \delta v^{\rm N}(0,a)
 =\frac{y_0^{\rm N}(a)}{\gamma}
       \left(\gamma a+\log y_0^{\rm N}(a)-1\right).
\end{align}
The non-positive value in \eqref{eq:inactive-lower-value} forces
$y_0^{\rm N}(a)\to0$ as $a\to\infty$.
If $a_2<\kappa/\gamma$, the right-hand side of
\eqref{eq:positive-endpoint-integrated} tends to $-\infty$, contradicting
$a\to\infty$.  If $a_2>\kappa/\gamma$, substituting
\eqref{eq:positive-endpoint-integrated} into
\eqref{eq:inactive-lower-value} gives
\[
 \delta v^{\rm N}(0,a)
 =\left(a_2-\frac{\kappa}{\gamma}\right)
       y_0^{\rm N}(a)^{1-1/\kappa}
  +\frac{\kappa-1}{\gamma}y_0^{\rm N}(a)\longrightarrow+\infty, ~~~~~ \text{as}~~ a \to \infty,\]
because $\kappa<1$ by Lemma \ref{lem:boundary-ratio}.
This contradicts the non-positive value, therefore, $a_2=\kappa/\gamma$, and
\eqref{eq:positive-endpoint-integrated} yields
\begin{align}\label{eq:positive-endpoint-solution}
 y_0^{\rm N}(a)=e^{\kappa-\gamma a},\qquad
 y_b^{\rm N}(a)=\frac{e^{\kappa-\gamma a}}{z_{0b}},
 \qquad a\geq a_2.
\end{align}

Substituting the solution
\eqref{eq:dual-general-r-low} into the two contact conditions
\eqref{eq:positive-rate-contact} gives
\begin{align}
 &\nu_+C_1^{\rm N}(a)y_0^{\rm N}(a)^{\nu_+-1}
 +\nu_-C_2^{\rm N}(a)y_0^{\rm N}(a)^{\nu_--1}
 =-\frac{\log y_0^{\rm N}(a)+\gamma a}{\gamma r}
   -\frac{m+\delta-r}{\gamma r^2}, \\
& \nu_+(\nu_+-1)C_1^{\rm N}(a)y_0^{\rm N}(a)^{\nu_+-2}
 +\nu_-(\nu_--1)C_2^{\rm N}(a)y_0^{\rm N}(a)^{\nu_--2}
 =-\frac1{\gamma ry_0^{\rm N}(a)}.
\end{align}
By solving the above two equations, we obtain
\begin{align}
 C_1^{\rm N}(a)
 =-\frac{(1-r\theta)(1-\nu_-)z_{0b}^{\nu_+-1}}
          {\gamma r(\nu_+-1)(\nu_+-\nu_-)}
       y_0^{\rm N}(a)^{1-\nu_+},\quad
 C_2^{\rm N}(a)
 =\frac{(1-r\theta)(\nu_+-1)z_{0b}^{\nu_--1}}
          {\gamma r(1-\nu_-)(\nu_+-\nu_-)}
       y_0^{\rm N}(a)^{1-\nu_-},\\
       \label{eq:inactive-positive-coefficientC1}
\end{align}
with $ C_1^{\rm N}(a)<0$ and $ C_2^{\rm N}(a)>0$. Differentiation with $(y_0^{\rm N})'(a)=-\gamma y_0^{\rm N}(a)$
recovers \eqref{eq:coefficient-derivatives13}--\eqref{eq:coefficient-derivatives24}, verifying the two
annuitization conditions.

Finally, substituting the coefficients $C_1^{\rm N}(a),  C_2^{\rm N}(a)$, and
$y_b^{\rm N}(a)=y_0^{\rm N}(a)/z_{0b}$ into
$b^{\rm N}(a)=-\widehat v_y^{\rm N}(y_b^{\rm N}(a),a)$, using $\log y_0^{\rm N}(a)=\kappa-\gamma a$ and  \eqref{eq:general-roots-v+}, and further simplifying the terms, we obtain
\begin{align}\label{eq:inactive-boundary-bN}
 b^{\rm N}(a)
 =\frac{1-\kappa+\log z_{0b}-\theta(m+\delta)}{\gamma r}>0.
\end{align}
Positivity of $b^{\rm N}(a)$ will be proved in the following
Lemma~\ref{lem:positive-inactive-signs}. Let $x_c^{\rm N}(a)$ denote the zero-consumption threshold. Since in this case, at $a=a_2$
consumption is zero at $x=0$ and positive at every positive wealth
level in the continuation region, thus $x_c^{\rm N}(a) =0$ for $a \ge a_2$.
Moreover, at zero wealth, we have
\begin{align}\label{eq:inactive-c-at-zero}
 c^{\rm N}(0,a)=-\frac1\gamma\log y_0^{\rm N}(a)
 =a-\frac{\kappa}{\gamma}<a,
 \qquad \pi^{\rm N}(0,a)=0.
\end{align}
Thus, the lower
state-constraint conditions in \eqref{eq:lower_state_constraint} are satisfied in this case. 
\subsubsection{Mixed consumption: $(\kappa/{\gamma})z_{0b}^{-1/\kappa}<a<\kappa/\gamma$ }
If the annuity income $0 \le a < a_2$, the dual endpoint $y^{\rm N}_0(a)$
would exceed one and the unconstrained consumption rate at zero wealth
would become negative.  We therefore use both the dual value functions given in \eqref{eq:dual-general-r-low} and \eqref{eq:dual-general-r-high}. Assume that 
$0<y_b^{\rm N}(a)<1<y_0^{\rm N}(a)$. The consumption boundary $x_c^{\rm N}(a)$ is determined by $y=1$.

In this case, we have $\widetilde U'(y_0^{\rm N}(a))=0$, so
\eqref{eq:uncapped-endpoint-equation} reduces to
\begin{align}\label{eq:y0_derivaive_Ninjection}
 \frac{(y_0^{\rm N})'(a)}{y_0^{\rm N}(a)}=-\frac{\kappa}{a}.
\end{align}
Integrating from $a$ to $a_2$ and imposing
$y_0^{\rm N}(a_2)=1$ give
$- \log y_0^{\rm N}(a)= \kappa\log\frac{a}{a_2}$, and thus,
\begin{align}\label{eq:mixed-endpoint-solution}
 y_0^{\rm N}(a)=\left(\frac{\kappa}{\gamma a}\right)^\kappa,
 \qquad y_b^{\rm N}(a)=\frac1{z_{0b}}
       \left(\frac{\kappa}{\gamma a}\right)^\kappa.
\end{align}
At the zero-wealth boundary, substituting
\eqref{eq:dual-general-r-high} into \eqref{eq:positive-rate-contact}
gives
\begin{align}
 &\nu_+C_3^{\rm N}(a)y_0^{\rm N}(a)^{\nu_+-1}
 +\nu_-C_4^{\rm N}(a)y_0^{\rm N}(a)^{\nu_--1}
 =-\frac ar,  \label{eq:positive-rate-C341}\\
 &\nu_+(\nu_+-1)C_3^{\rm N}(a)y_0^{\rm N}(a)^{\nu_+-2}
 +\nu_-(\nu_--1)C_4^{\rm N}(a)y_0^{\rm N}(a)^{\nu_--2}
 =0. \label{eq:positive-rate-C342}
\end{align}
Solving \eqref{eq:positive-rate-C341} and \eqref{eq:positive-rate-C342}, we obtain
\begin{align}\label{eq:inactive-zero-coefficient-C3}
 C_3^{\rm N}(a)
 =-\frac{(1-\nu_-)a}{r\nu_+(\nu_+-\nu_-)}
       y_0^{\rm N}(a)^{1-\nu_+},\quad
 C_4^{\rm N}(a)
 =-\frac{(\nu_+-1)a}{r\nu_-(\nu_+-\nu_-)}
       y_0^{\rm N}(a)^{1-\nu_-}. 
\end{align}
The coefficients on the positive-consumption side then follow from
the matching relations \eqref{eq:general-r-matching} at $y=1$:
\begin{align}\label{eq:inactive-mixed-coefficients}
 C_1^{\rm N}(a)&=C_3^{\rm N}(a)+\Delta_+,\qquad
 C_2^{\rm N}(a)=C_4^{\rm N}(a)+\Delta_-,
\end{align}
in which $\Delta_+$ and $\Delta_-$ are given in \eqref{eq:general-r-matching-data+1}. The dual equation \eqref{eq:dual-ode} and continuity of
$\widetilde U,\widetilde U'$ in \eqref{eq:utility-conjugate}
then give matching of the second $y$-derivatives at $y=1$.
Differentiating the solved coefficients and using
\eqref{eq:kappa-root-identity1} recovers
\eqref{eq:coefficient-derivatives13}--\eqref{eq:coefficient-derivatives24}; hence the annuitization conditions  \eqref{eq:dual-boundaries-annuitization-condition1} 
hold.

Differentiating \eqref{eq:dual-general-r-low} and
\eqref{eq:dual-general-r-high}, we get the wealth equations
on $y_b^{\rm N}(a)\leq y\leq1$ and
$1\leq y\leq y_0^{\rm N}(a)$, respectively:
\begin{align}\label{eq:wealth-map-positive}
 x={}&-\nu_+C_1^{\rm N}(a)y^{\nu_+-1}
      -\nu_-C_2^{\rm N}(a)y^{\nu_--1}
      -\frac ar-\frac{\log y}{\gamma r}
      -\frac{m+\delta-r}{\gamma r^2},\\
 x={}&-\nu_+C_3^{\rm N}(a)y^{\nu_+-1}
      -\nu_-C_4^{\rm N}(a)y^{\nu_--1}
      -\frac ar.
       \label{eq:wealth-map-zero}
\end{align}
Furthermore, setting $y=y_b^{\rm N}(a)=  y^{\rm N}_0(a)/z_{0b} $ in
\eqref{eq:wealth-map-positive} and using \eqref{eq:general-r-matching-data+}, we deduce the annuitization boundary $b^{\rm N}(a)$:
\begin{align}
 b^{\rm N}(a)={}&\frac ar\left[
 \frac{1-\nu_-}{\nu_+-\nu_-}z_{0b}^{1-\nu_+}
 +\frac{\nu_+-1}{\nu_+-\nu_-}z_{0b}^{1-\nu_-}-1
 \right]-\frac{\log(y_b^{\rm N}(a))}{\gamma r}
  -\frac{m+\delta-r}{\gamma r^2}\\
 &+\frac1{\gamma r(\nu_+-\nu_-)}
 \left[
 \frac{1-\nu_-}{\nu_+-1}
       \left({y_b^{\rm N}(a)}\right)^{\nu_+-1}
 -\frac{\nu_+-1}{1-\nu_-}
       \left({y_b^{\rm N}(a)}\right)^{\nu_--1}
 \right].\label{eq:inactive-mixed-annuity-boundary}
\end{align}
Setting $y=1$ in \eqref{eq:wealth-map-zero}, we obtain the
consumption boundary $x_c^{\rm N}(a)$:
\begin{align}\label{eq:inactive-mixed-consum-boundary}
 x_c^{\rm N}(a)
 =\frac ar\left[
 \frac{1-\nu_-}{\nu_+-\nu_-}y_0^{\rm N}(a)^{1-\nu_+}
 +\frac{\nu_+-1}{\nu_+-\nu_-}y_0^{\rm N}(a)^{1-\nu_-}-1
 \right].
\end{align}
In particular, at $a_2 = \kappa/\gamma$, $y_0^{\rm N}(a_2)=1$ and $x_c^{\rm N}(a_2)=0$.
Assume that at the lower annuity income threshold $a = a_1$, we have $y_b^{\rm N}(a_1)=1$. Substituting this condition into
\eqref{eq:mixed-endpoint-solution}, we obtain $a_1=\kappa z_{0b}^{-1/\kappa}/\gamma<a_2$, which yields that $x_c^{\rm N}(a_1)=b^{\rm N}(a_1)$.
\subsubsection{Zero consumption: $0<a \le (\kappa/{\gamma})z_{0b}^{-1/\kappa}$ }
For $0<a\leq a_1$, we conjecture that the optimal
consumption is identically zero over the
entire continuation region:
$1\leq y_b^{\rm N}(a)<y_0^{\rm N}(a)$, with
$y_b^{\rm N}(a_1)=1$.
Using $\widetilde U(y)=-1/\gamma$ from \eqref{eq:utility-conjugate}
in \eqref{eq:dual-ode}, the continuation equation becomes
\[
 my^2\widehat v_{yy}^{\rm N}
 +(\delta-r)y\widehat v_y^{\rm N}
 -\delta\widehat v^{\rm N}+ay-\frac1\gamma=0,
\]
and the dual value function is given in \eqref{eq:dual-general-r-high}.
In view that $\widetilde U'(y_0^{\rm N}(a))=0$, the endpoint $y_0^{\rm N}(a)$ satisfies the differential equation \eqref{eq:y0_derivaive_Ninjection} for $0<a\leq a_1$. Integrating this equation subject to $y_0^{\rm N}(a_1)=z_{0b}$ yields
\[
 y_0^{\rm N}(a)
 =z_{0b}\left(\frac{a_1}{a}\right)^\kappa
 =\left(\frac{\kappa}{\gamma a}\right)^\kappa,
\]
where the last equality follows from $a_1=\kappa z_{0b}^{-1/\kappa}/\gamma$. As in the mixed case,
the contact equations \eqref{eq:positive-rate-C341}--\eqref{eq:positive-rate-C342}
remain unchanged. Substituting $y_0^{\rm N}(a)$ into these equations
gives the coefficients $C_3^{\rm N}(a)$ and $C_4^{\rm N}(a)$ in \eqref{eq:inactive-zero-coefficient-C3}.

Furthermore,
substituting these coefficients and
$y=y_b^{\rm N}(a)=y_0^{\rm N}(a)/z_{0b}$
into \eqref{eq:wealth-map-zero} gives
\begin{align}\label{eq:inactive-zero-boundaries}
 b^{\rm N}(a)=x_c^{\rm N}(a)
 =\frac ar\left[
 \frac{1-\nu_-}{\nu_+-\nu_-}z_{0b}^{1-\nu_+}
 +\frac{\nu_+-1}{\nu_+-\nu_-}z_{0b}^{1-\nu_-}-1
 \right].
\end{align}
\begin{remark}
Combining the exponential endpoint \eqref{eq:positive-endpoint-solution}
with the power-law endpoint \eqref{eq:mixed-endpoint-solution}
gives
\begin{align}\label{eq:y0-positive}
 y_0^{\rm N}(a)=
 \begin{cases}
 \left(\dfrac{\kappa}{\gamma a}\right)^\kappa,&0<a< a_2,\\[0.7em]
 e^{\kappa-\gamma a},&a\geq a_2,
 \end{cases}
\end{align}
The two expressions in \eqref{eq:y0-positive} have the same
value and first derivative at $a=a_2$.
As $a\downarrow0$, $y_0^{\rm N}(a)\to\infty$ and
$b^{\rm N}(a)\to0$, yielding the degenerate no-injection policy at
$(0,0)$.  The funded problem in Section \ref{sec:injection-active} below instead has finite dual endpoints
at $a=0$.
\end{remark}
\paragraph{Convexity and the state constraint.}
To justify the inverse transform, we first establish strict dual
convexity. Define $h^{\rm N}(y,a)=y\widehat v_{yy}^{\rm N}(y,a),
 $ for $y_b^{\rm N}(a)\leq y\leq y_0^{\rm N}(a)$. This is also the factor entering the risky-investment rule.

\begin{lemma}[Properties of the inactive dual solution]\label{lem:positive-inactive-signs}
For every $a>0$, the solution constructed above satisfies $h^{\rm N}(y_0^{\rm N}(a),a)=0$, $h_y^{\rm N}(y,a)<0$, and 
$h^{\rm N}(y,a)>0$ for $y_b^{\rm N}(a)\leq y<y_0^{\rm N}(a)$. The dual solution $\widehat v^{\rm N}(y,a)$ is strictly convex
and decreases in $y$ in the interior and satisfies $0\leq\widehat v_{ay}^{\rm N}\leq\theta$ and
 $0<\widehat v_a^{\rm N}\leq\theta y$. Moreover, $b^{\rm N}(a)>0$, $(b^{\rm N})'(a)>-\theta$, and $\pi^{\rm N}(0,a)=0$, $0\leq c^{\rm N}(0,a)<a$.
\end{lemma}

\subsection{Solving the injection-active problem}\label{sec:injection-active}
We now allow capital injection and use the superscript
$\mathrm{I}$ for the resulting solution.  The no-injection endpoint $y_0^{\rm N}(a)$ in
\eqref{eq:y0-positive} continuously decreases from infinity to zero.
Thus, there is a unique income level $a_\beta$ at which $y_0^{\rm N}(a_\beta)=\beta$. Solving this condition separately for each branch of \eqref{eq:y0-positive} yields
\begin{align}\label{eq:abeta-positive}
 a_\beta=
 \begin{cases}
 \dfrac{\kappa}{\gamma\beta^{1/\kappa}},&\beta\geq1,\\[0.8em]
 \dfrac{\kappa-\log\beta}{\gamma},&0<\beta<1.
 \end{cases}
\end{align}
For $a\geq a_\beta$, then $y_0^{\rm N}(a) \le \beta$ and the no-injection solution satisfies the marginal constraint without reflection and therefore remains unchanged. For $a<a_\beta$, the boundary condition $v_x(0,a)=\beta$ binds, and we replace the super-contact condition \eqref{eq:dual-state-constraint} with $y_0^{\rm I}(a)=\beta$. We then solve the boundary equations subject to this constraint and examine the three consumption cases.

\paragraph{The endpoints and coefficients.}
The ratio equation \eqref{eq:z0b-general} holds under either the
super-contact condition or the binding marginal cap.  Consequently,
\begin{align}\label{eq:active-endpoints}
 y_0^{\rm I}(a)=\beta,\qquad
 y_b^{\rm I}(a)=\frac{\beta}{z_{0b}},\qquad 0\leq a<a_\beta.
\end{align}
Substituting \eqref{eq:active-endpoints} into
\eqref{eq:coefficient-derivatives13}--\eqref{eq:coefficient-derivatives24}
makes the derivatives of the coefficients constant in $a$.  Integrating from
$a_\beta$, with $C_j^{\rm I}(a_\beta)=C_j^{\rm N}(a_\beta)$ yields
\begin{align}
 C_j^{\rm I}(a)
 &=C_j^{\rm N}(a_\beta)
 -\frac{1-r\theta}{r}\frac{1-\nu_-}{\nu_+-\nu_-}
       \left(\frac{\beta}{z_{0b}}\right)^{1-\nu_+}(a-a_\beta),
       &&j=1,3,\label{eq:active-coeff-positive13}\\
 C_j^{\rm I}(a)
 &=C_j^{\rm N}(a_\beta)
 -\frac{1-r\theta}{r}\frac{\nu_+-1}{\nu_+-\nu_-}
       \left(\frac{\beta}{z_{0b}}\right)^{1-\nu_-}(a-a_\beta),
       &&j=2,4. \label{eq:active-coeff-positive24}
\end{align}
At $a_\beta$, the no-injection coefficients are given by
\eqref{eq:inactive-positive-coefficientC1} when $a_\beta\geq a_2$,
or by \eqref{eq:inactive-zero-coefficient-C3} and \eqref{eq:inactive-mixed-coefficients} when $a_\beta<a_2$.
Equations \eqref{eq:active-coeff-positive13} and  \eqref{eq:active-coeff-positive24} give
$(C_1^{\rm I})'=(C_3^{\rm I})'$ and $(C_2^{\rm I})'=(C_4^{\rm I})'$.
The differences between the two pairs of coefficients therefore remain
$\Delta_+$ and $\Delta_-$ given in \eqref{eq:general-r-matching-data+1}, preserving value and slope matching
at $y=1$.

Combining the active and inactive ranges, we get
\begin{align}\label{eq:positive-rate-capped-endpoints}
 y_0^{\rm I}(a)=\min\{\beta,y_0^{\rm N}(a)\},\qquad
 y_b^{\rm I}(a)=\frac{y_0^{\rm I}(a)}{z_{0b}}.
\end{align}
At $a=0$, the minimum in \eqref{eq:positive-rate-capped-endpoints}
is understood through $y_0^{\rm N}(a)\to\infty$ as $a\downarrow0$.
For $a\geq a_\beta$, set $C_j^{\rm I}(a)=C_j^{\rm N}(a), j=1,2,3,4$.
The coefficient functions coincide at $a_\beta$ by construction.
Their first derivatives coincide by
\eqref{eq:coefficient-derivatives13}--\eqref{eq:coefficient-derivatives24},
because both sides have the same endpoint
$y_b(a_\beta)=\beta/z_{0b}$.

\paragraph{The boundary conditions and the second derivative of the dual value function.}
Substituting \eqref{eq:active-coeff-positive13} and \eqref{eq:active-coeff-positive24} into the dual branches
\eqref{eq:dual-general-r-low}--\eqref{eq:dual-general-r-high}, after a simple calculation, we obtain
\begin{align}\label{eq:active-dual-continuation}
 \widehat v^{\rm I}(y,a)
 =\widehat v^{\rm N}(y,a_\beta)
 +(a-a_\beta)\widehat v_a^{\rm N}(y,a_\beta),
 \qquad \frac{\beta}{z_{0b}}\leq y\leq\beta,\quad 0\leq a<a_\beta.
\end{align}
At $(y,a)=(\beta,a_\beta)$, from \eqref{eq:positive-rate-contact} and \eqref{eq:deriv_hv_ay} it follows that
$\widehat v_y^{\rm N}(\beta,a_\beta)
=\widehat v_{ay}^{\rm N}(\beta,a_\beta)=0$.  Differentiating
\eqref{eq:active-dual-continuation} in $y$  verifies
$\widehat v_y^{\rm I}(\beta,a)=0$.
The derivatives of the coefficients in
\eqref{eq:coefficient-derivatives13}--\eqref{eq:coefficient-derivatives24}
satisfy the annuitization condition
\eqref{eq:dual-boundaries-annuitization-condition1}.
Differentiating \eqref{eq:active-dual-continuation} twice in $y$ leads to
\begin{align}\label{eq:positive-rate-curvature-comparison}
 \widehat v_{yy}^{\rm I}(y,a)
 =\widehat v_{yy}^{\rm N}(y,a_\beta)
 +(a-a_\beta)\widehat v_{ayy}^{\rm N}(y,a_\beta)>0.
\end{align}
Indeed, the first term is non-negative by
Lemma~\ref{lem:positive-inactive-signs}, whereas
$\widehat v_{ayy}^{\rm N}=g_{yy}<0$ by
\eqref{eq:general-r-g-explicit}, and $a-a_\beta<0$.
Thus, \eqref{eq:positive-rate-curvature-comparison} establishes the strict
convexity of $\widehat{v}^I$ on ${\beta}/{z_{0b}}\leq y\leq\beta,\quad 0\leq a<a_\beta$.

At $y=\beta$, inserting \eqref{eq:positive-rate-contact} and
\eqref{eq:endpoint-mixed-curvature} into
\eqref{eq:positive-rate-curvature-comparison} gives
$\widehat v_{yy}^{\rm I}(\beta,a)=\kappa(a_\beta-a)/(m\beta)$.
Using the Legendre identities \eqref{eq:dual-relations} in the
regular feedback control \eqref{eq:feedback-regular}, we obtain
\begin{align}\label{eq:active-zero-wealth-controls}
 \pi^{\rm I}(0,a)=\frac{\mu-r}{\sigma^2}\frac{\kappa}{m}(a_\beta-a)>0,
 \qquad c^{\rm I}(0,a)=\frac{(-\log\beta)^+}{\gamma}.
\end{align}
The investment rule in \eqref{eq:active-zero-wealth-controls} tends
to zero as $a\uparrow a_\beta$, where lower reflection is no longer
needed.  At the annuitization endpoint, differentiate
\eqref{eq:active-dual-continuation} in $y$ and use
$\widehat v_{ay}^{\rm N}(\beta/z_{0b},a_\beta)=\theta$
from \eqref{eq:dual-boundaries-annuitization-condition1}.
The wealth relation $b^{\rm I}(a)=-\widehat v_y^{\rm I}(\beta/z_{0b},a)$
then yields
\begin{align}\label{eq:active-annuity-boundary}
 b^{\rm I}(a)=b^{\rm N}(a_\beta)+\theta(a_\beta-a),
 \qquad 0\leq a<a_\beta.
\end{align}
We now evaluate \eqref{eq:active-annuity-boundary} and the consumption
threshold on each branch.

\subsubsection{Positive consumption}
Suppose $0<\beta\leq1$.  By \eqref{eq:active-endpoints}, the whole
dual interval lies on the positive-consumption branch, and
\eqref{eq:abeta-positive} gives $a_\beta\geq a_2$.
The dual value function $\widehat v_+^{\rm I}$ is defined by \eqref{eq:dual-general-r-low} where $C_1^{\rm I}(a)$ and $C_2^{\rm I}(a)$ are given in \eqref{eq:active-coeff-positive13} and \eqref{eq:active-coeff-positive24}, respectively. Moreover,
substituting \eqref{eq:inactive-boundary-bN} into
\eqref{eq:active-annuity-boundary}, we obtain
\begin{align}\label{eq:active-positive-boundaries}
 b^{\rm I}(a)
 =\frac{1-\kappa+\log z_{0b}-\theta(m+\delta)}{\gamma r}
   +\theta(a_\beta-a),\qquad x_c^{\rm I}(a)=0.
\end{align}
The purchase boundary in \eqref{eq:active-positive-boundaries}
is linear in $a$, and consumption is positive at every
positive-wealth point of the continuation region.
If $\beta<1$, consumption is positive even at zero wealth;
capital injections may then finance consumption as well as offset investment
losses.  At $\beta=1$, consumption is zero at $x=0$ and positive
for every $x>0$ in the continuation region.

\subsubsection{Mixed consumption}
Suppose $1<\beta<z_{0b}$.  Then $\beta/z_{0b}<1<\beta$ and
$a_1<a_\beta<a_2$.  The dual value function equals $\widehat v_0^{\rm I}$ given in \eqref{eq:dual-general-r-high}  on $[1,\beta]$ and $\widehat v_+^{\rm I}$ given in \eqref{eq:dual-general-r-low} on $[\beta/z_{0b},1]$ with
the coefficients given in \eqref{eq:active-coeff-positive13} and \eqref{eq:active-coeff-positive24}, respectively.

The consumption boundary is the wealth corresponding to $y=1$.
Differentiating \eqref{eq:active-dual-continuation} in $y$ and
using $x=-\widehat v_y(y,a)$ at $y=1$ gives
$x_c^{\rm I}(a)
 =x_c^{\rm N}(a_\beta)
 +(a_\beta-a)\widehat v_{ay}^{\rm N}(1,a_\beta)$. Substituting \eqref{eq:annuity-mixed-derivative}, with
$y_b^{\rm N}(a_\beta)=\beta/z_{0b}$ therefore yields
\begin{align}
 x_c^{\rm I}(a)
 ={}&x_c^{\rm N}(a_\beta)+\frac{a_\beta-a}{r}\\
 &-\frac{(a_\beta-a)(1-r\theta)}{r}
 \left[
 \frac{\nu_+(1-\nu_-)}{\nu_+-\nu_-}
       \left(\frac{z_{0b}}{\beta}\right)^{\nu_+-1}
 +\frac{\nu_-(\nu_+-1)}{\nu_+-\nu_-}
       \left(\frac{z_{0b}}{\beta}\right)^{\nu_--1}
 \right].\label{eq:active-mixed-boundaries}
\end{align}
Here, $x_c^{\rm N}(a_\beta)$ is obtained by evaluating
\eqref{eq:inactive-mixed-consum-boundary} at $a=a_\beta$.
The purchase boundary follows from \eqref{eq:active-annuity-boundary}
and \eqref{eq:inactive-mixed-annuity-boundary}.
Strict convexity of the dual value function implies that the threshold in
\eqref{eq:active-mixed-boundaries} satisfies
$0<x_c^{\rm I}(a)<b^{\rm I}(a)$.
Within the continuation region, consumption is zero for
$0\leq x\leq x_c^{\rm I}(a)$ and positive for
$x_c^{\rm I}(a)<x<b^{\rm I}(a)$.
\subsubsection{Zero consumption}
Suppose $\beta\geq z_{0b}$.  Then $\beta/z_{0b}\geq1$ and
$a_\beta\leq a_1$, so the dual value function equals $\widehat v_0^{\rm I}$ with $C_3^{\rm I}(a)$ and $C_4^{\rm I}(a)$ given in \eqref{eq:active-coeff-positive13} and \eqref{eq:active-coeff-positive24}, respectively. The boundary values are specified by \eqref{eq:inactive-zero-coefficient-C3} at $a_\beta$. Substitution of \eqref{eq:inactive-zero-boundaries} into
\eqref{eq:active-annuity-boundary} gives
\begin{align}\label{eq:active-zero-boundaries}
\begin{split}
 x_c^{\rm I}(a)=b^{\rm I}(a)
 ={}&\frac{a_\beta}{r}
 \left[
 \frac{1-\nu_-}{\nu_+-\nu_-}z_{0b}^{1-\nu_+}
 +\frac{\nu_+-1}{\nu_+-\nu_-}z_{0b}^{1-\nu_-}-1
 \right]
 +\theta(a_\beta-a).
\end{split}
\end{align}
Thus, in this case, consumption is zero throughout the continuation region: the injection supports the risky-investment strategy when wealth is zero, but it does not finance current consumption.

\subsection{The primal solution and the six policy cases}\label{sec:six-regimes}
In this subsection, we characterize the value function and regular controls from
the dual solution, summarize the six cases, and extend the
candidate to the annuitization region.
Unless a comparison requires the superscript $\mathrm{I}$, we write
$\widehat v$, $C_j(a)$, $y_0(a)$, $y_b(a)$, $b(a)$, and $x_c(a)$
for the problem with capital injection.  In particular, the dual endpoints are
those given in \eqref{eq:positive-rate-capped-endpoints};
$y_0^{\rm N}(a)$ continues to denote the no-injection endpoint.

\subsubsection{The inverse transform and regular controls}
The policy boundaries are obtained from the relationship between wealth and the dual variable as follows:
\begin{align}\label{eq:positive-boundary-recovery-ba}
 b(a)=-\widehat v_y(y_b(a),a),\quad
 x_c(a)=
 \begin{cases}
 b(a),&y_b(a)\geq1,\\
 -\widehat v_y(1,a),&y_b(a)<1<y_0(a),\\
 0,&y_0(a)\leq1.  
 \end{cases}
\end{align}
For $a\geq a_\beta$, in the positive-consumption case, the annuity purchase boundary equals
\eqref{eq:inactive-boundary-bN}, with $x_c(a)=0$;
in the mixed case, the purchase and consumption boundaries are
\eqref{eq:inactive-mixed-annuity-boundary} and
\eqref{eq:inactive-mixed-consum-boundary}, respectively;
in the zero-consumption case, both are given by
\eqref{eq:inactive-zero-boundaries}.
For $a<a_\beta$, the corresponding formulas are
\eqref{eq:active-positive-boundaries} in the positive case,
\eqref{eq:active-annuity-boundary} and
\eqref{eq:active-mixed-boundaries} in the mixed case,
and \eqref{eq:active-zero-boundaries} in the zero case.

For $0\leq x\leq b(a)$, define $y^*(x,a)$ as the solution of
\begin{align}\label{eq:positive-dual-inversion}
 x=-\widehat v_y(y^*(x,a),a),\qquad
 y_b(a)\leq y^*(x,a)\leq y_0(a).
\end{align}
In particular, for $y^*(x,a) \le 1$ and $y^*(x,a) \ge 1$, we have
\begin{align}
 x={}&-\nu_+C_1(a)(y^*(x,a))^{\nu_+-1}
      -\nu_-C_2(a)(y^*(x,a))^{\nu_--1}
      -\frac ar-\frac{\log y^*(x,a)}{\gamma r}
      -\frac{m+\delta-r}{\gamma r^2},\\
 x={}&-\nu_+C_3(a)(y^*(x,a))^{\nu_+-1}
      -\nu_-C_4(a)(y^*(x,a))^{\nu_--1}
      -\frac ar,
 \end{align}
in which $C_j, j=1,2,3,4.$ are given in \eqref{eq:active-coeff-positive13} and \eqref{eq:active-coeff-positive24}.
Lemma~\ref{lem:positive-inactive-signs} and
\eqref{eq:positive-rate-curvature-comparison} show that
as $y$ increases from $y_b(a)$ to $y_0(a)$, $-\widehat{v}_y(y,a)$ is strictly decreasing from $b(a)$ to zero.
Hence, \eqref{eq:positive-dual-inversion} has exactly one solution.

Define the primal candidate function by
\begin{align}\label{eq:inverse-transform}
 v(x,a)=\inf_{y_b(a)\leq y\leq y_0(a)}
       \{\widehat v(y,a)+xy\}.
\end{align}
The minimizer in \eqref{eq:inverse-transform} is $y^*(x,a)$.
In \eqref{eq:primal-branches-explicit} and
\eqref{eq:investment-branches-explicit}, we write $y=y^*(x,a)$.
Substituting $x=-\widehat v_y(y,a)$
into \eqref{eq:dual-general-r-low}--\eqref{eq:dual-general-r-high} and canceling the linear terms,
we obtain
\begin{align}\label{eq:primal-branches-explicit}
 v(x,a)=
 \begin{cases}
 (1-\nu_+)C_1(a)y^{\nu_+}
 +(1-\nu_-)C_2(a)y^{\nu_-}-\dfrac{y}{\gamma r},
       &y\leq1,\\[0.7em]
 (1-\nu_+)C_3(a)y^{\nu_+}
 +(1-\nu_-)C_4(a)y^{\nu_-}-\dfrac1{\gamma\delta},
       &y\geq1.
 \end{cases}
\end{align}
The branches in \eqref{eq:primal-branches-explicit} agree at $y=1$
by \eqref{eq:general-r-matching}.
Substituting the Legendre identities \eqref{eq:dual-relations}
into \eqref{eq:feedback-regular}, we obtain the regular controls
\begin{align}\label{eq:dual-feedback-detailed}
 c^*(x,a)=\frac{(-\log y^*(x,a))^+}{\gamma},\qquad
 \pi^*(x,a)=\frac{\mu-r}{\sigma^2}
      y^*(x,a)\widehat v_{yy}(y^*(x,a),a).
\end{align}
Moreover, differentiating \eqref{eq:dual-general-r-low}--\eqref{eq:dual-general-r-high} twice with respect to $y$, we deduce the explicit expression of the investment rule:
\begin{align}\label{eq:investment-branches-explicit}
 \frac{\sigma^2}{\mu-r}\pi^*(x,a)=
 \begin{cases}
 \nu_+(\nu_+-1)C_1(a)y^{\nu_+-1}
 +\nu_-(\nu_--1)C_2(a)y^{\nu_--1}+\dfrac1{\gamma r},
       &y\leq1,\\[0.7em]
 \nu_+(\nu_+-1)C_3(a)y^{\nu_+-1}
 +\nu_-(\nu_--1)C_4(a)y^{\nu_--1},
       &y\geq1.
 \end{cases}
\end{align}

\begin{proposition}[Inversion of the dual problem]\label{prop:dual-primal}
The candidate in \eqref{eq:inverse-transform} is increasing and strictly
concave in $x$ on $0<x<b(a)$, and increasing in $a$.
It solves the continuation equation \eqref{eq:primal-continuation} and satisfies the annuity-purchase inequality
\[
v_a(x,a)-\theta v_x(x,a)\leq 0,
\qquad 0\leq x\leq b(a).
\]
Moreover, it satisfies the zero-wealth boundary condition
$v_x(0,a)=y_0(a)\leq\beta$ and the annuitization contact conditions
\begin{align}\label{eq:primal-smooth-fit}
 v_a(b(a),a)=\theta v_x(b(a),a),\qquad
 v_{ax}(b(a),a)=\theta v_{xx}(b(a),a).
\end{align}
\end{proposition}

\subsubsection{Six policy cases}
The threshold $a_\beta$ defined in \eqref{eq:abeta-positive}
separates the regions where capital injection is active
($a<a_\beta$) and inactive ($a\geq a_\beta$).
The consumption regime is determined by the position of $1$
relative to $y_b(a)$ and $y_0(a)$.
Combining the two injection regimes with the three consumption
regimes gives six cases, summarized in Table~\ref{tab:six-regimes}.

When capital injection is inactive, the coefficients on the relevant
consumption branches are given by
\eqref{eq:inactive-positive-coefficientC1},
\eqref{eq:inactive-zero-coefficient-C3},
and \eqref{eq:inactive-mixed-coefficients}.
When capital injection is active, the coefficients are given by
\eqref{eq:active-coeff-positive13} and \eqref{eq:active-coeff-positive24}.
Substituting the corresponding coefficients into
\eqref{eq:primal-branches-explicit} and
\eqref{eq:dual-feedback-detailed} yields the value function
and regular controls, respectively.
In Figure~\ref{fig:six-slices}, we illustrate the zero-consumption,
positive-consumption, and immediate-annuitization regions along
the wealth axis for each of the six cases. Orange, blue, and green indicate zero consumption, positive
consumption, and immediate annuitization, respectively.
A filled red point marks reflection at zero wealth, where
$v_x(0,a)=\beta$; an open point marks a state constraint
satisfied without capital injection. In Figure~\ref{fig:policy-regions}, we plot the policy
regions and their boundaries in the $(a,x)$-plane when $a_\beta<a_1<a_2$ and $\beta>z_{0b}$.
 The vertical dashed line separates the
injection-active and inactive ranges.  Injection occurs only at $x=0$ in
the active range (red segment); the gray segment is maintained without
injection.  Orange and blue denote zero and positive consumption in the
continuation region, and green denotes immediate annuitization.  The solid
curve is $b(a)$, which is constant for $a\geq a_2$; the blue dashed curve is $x_c(a)$, which coincides with
$b(a)$ for $a\leq a_1$ and with zero for $a\geq a_2$. 

\begin{table}[htbp]
\centering
\small
\renewcommand{\arraystretch}{1.2}
\begin{tabular}{llll}
\hline
Injection & Consumption & Parameter range & Dual solution\\
\hline
Inactive & positive & $a\geq\max\{a_\beta,a_2\}$ & $\widehat v_+^{\rm N}$\\
Inactive & mixed & $a\geq a_\beta,\ a_1<a<a_2$
 & $\widehat v_+^{\rm N},\widehat v_0^{\rm N}$\\
Inactive & zero & $a_\beta\leq a\leq a_1$ & $\widehat v_0^{\rm N}$\\
Active & positive & $0\leq a<a_\beta,\ \beta\leq1$
 & $\widehat v_+^{\rm I}$\\
Active & mixed & $0\leq a<a_\beta,\ 1<\beta<z_{0b}$
 & $\widehat v_+^{\rm I},\widehat v_0^{\rm I}$\\
Active & zero & $0\leq a<a_\beta,\ \beta\geq z_{0b}$
 & $\widehat v_0^{\rm I}$\\
\hline
\end{tabular}
\caption{The six policy cases with two injection and three consumption
regimes.}
\label{tab:six-regimes}
\end{table}


\begingroup\color{black}
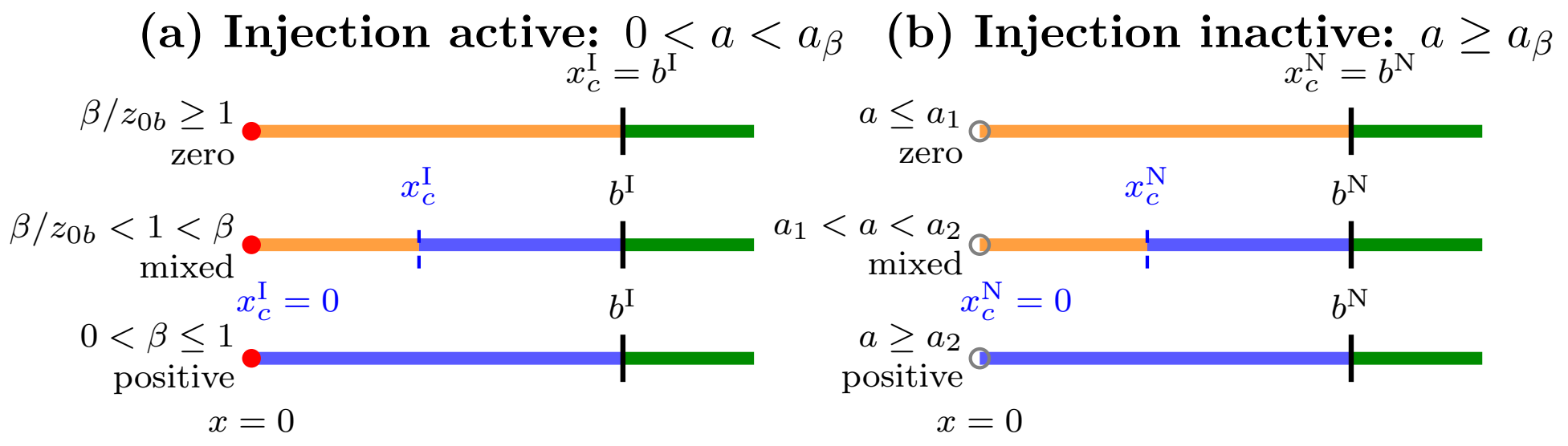
\begin{figure}[H]
\centering
\begin{tikzpicture}[x=1.20cm,y=0.78cm]
  \node[font=\bfseries] at (2.65,4.55) {(a) Injection active: $0<a<a_\beta$};
  \node[font=\bfseries] at (7.65,4.55) {(b) Injection inactive: $a\geq a_\beta$};
  \node[left,align=right,font=\scriptsize] at (1.0,3.55)
       {$\beta/z_{0b}\geq1$\\zero};
  \draw[line width=3pt,orange!75] (1.0,3.55)--(3.55,3.55);
  \draw[line width=3pt,green!55!black] (3.55,3.55)--(4.45,3.55);
  \fill[red] (1.0,3.55) circle (2.2pt);
  \draw[very thick] (3.55,3.30)--(3.55,3.80);
  \node[above,font=\scriptsize] at (3.55,3.80) {$x_c^{\rm I}=b^{\rm I}$};
  \node[left,align=right,font=\scriptsize] at (1.0,2.35)
       {$\beta/z_{0b}<1<\beta$\\mixed};
  \draw[line width=3pt,orange!75] (1.0,2.35)--(2.15,2.35);
  \draw[line width=3pt,blue!65] (2.15,2.35)--(3.55,2.35);
  \draw[line width=3pt,green!55!black] (3.55,2.35)--(4.45,2.35);
  \fill[red] (1.0,2.35) circle (2.2pt);
  \draw[thick,blue,dashed] (2.15,2.10)--(2.15,2.60);
  \draw[very thick] (3.55,2.10)--(3.55,2.60);
  \node[above,font=\scriptsize,blue] at (2.15,2.60) {$x_c^{\rm I}$};
  \node[above,font=\scriptsize] at (3.55,2.60) {$b^{\rm I}$};
  \node[left,align=right,font=\scriptsize] at (1.0,1.15)
       {$0<\beta\leq1$\\positive};
  \draw[line width=3pt,blue!65] (1.0,1.15)--(3.55,1.15);
  \draw[line width=3pt,green!55!black] (3.55,1.15)--(4.45,1.15);
  \fill[red] (1.0,1.15) circle (2.2pt);
  \draw[very thick] (3.55,0.90)--(3.55,1.40);
  \node[above,font=\scriptsize,blue] at (1.25,1.38) {$x_c^{\rm I}=0$};
  \node[above,font=\scriptsize] at (3.55,1.40) {$b^{\rm I}$};
  \node[left,align=right,font=\scriptsize] at (6.0,3.55)
       {$a\leq a_1$\\zero};
  \draw[line width=3pt,orange!75] (6.0,3.55)--(8.55,3.55);
  \draw[line width=3pt,green!55!black] (8.55,3.55)--(9.45,3.55);
  \draw[gray,thick] (6.0,3.55) circle (2.2pt);
  \draw[very thick] (8.55,3.30)--(8.55,3.80);
  \node[above,font=\scriptsize] at (8.55,3.80) {$x_c^{\rm N}=b^{\rm N}$};
  \node[left,align=right,font=\scriptsize] at (6.0,2.35)
       {$a_1<a<a_2$\\mixed};
  \draw[line width=3pt,orange!75] (6.0,2.35)--(7.15,2.35);
  \draw[line width=3pt,blue!65] (7.15,2.35)--(8.55,2.35);
  \draw[line width=3pt,green!55!black] (8.55,2.35)--(9.45,2.35);
  \draw[gray,thick] (6.0,2.35) circle (2.2pt);
  \draw[thick,blue,dashed] (7.15,2.10)--(7.15,2.60);
  \draw[very thick] (8.55,2.10)--(8.55,2.60);
  \node[above,font=\scriptsize,blue] at (7.15,2.60) {$x_c^{\rm N}$};
  \node[above,font=\scriptsize] at (8.55,2.60) {$b^{\rm N}$};
  \node[left,align=right,font=\scriptsize] at (6.0,1.15)
       {$a\geq a_2$\\positive};
  \draw[line width=3pt,blue!65] (6.0,1.15)--(8.55,1.15);
  \draw[line width=3pt,green!55!black] (8.55,1.15)--(9.45,1.15);
  \draw[gray,thick] (6.0,1.15) circle (2.2pt);
  \draw[very thick] (8.55,0.90)--(8.55,1.40);
  \node[above,font=\scriptsize,blue] at (6.25,1.38) {$x_c^{\rm N}=0$};
  \node[above,font=\scriptsize] at (8.55,1.40) {$b^{\rm N}$};
  \node[below,font=\scriptsize] at (1.0,0.78) {$x=0$};
  \node[below,font=\scriptsize] at (6.0,0.78) {$x=0$};
\end{tikzpicture}
\caption{The six policy cases with consumption and annuitization regions along the wealth axis.
}
\label{fig:six-slices}
\end{figure}

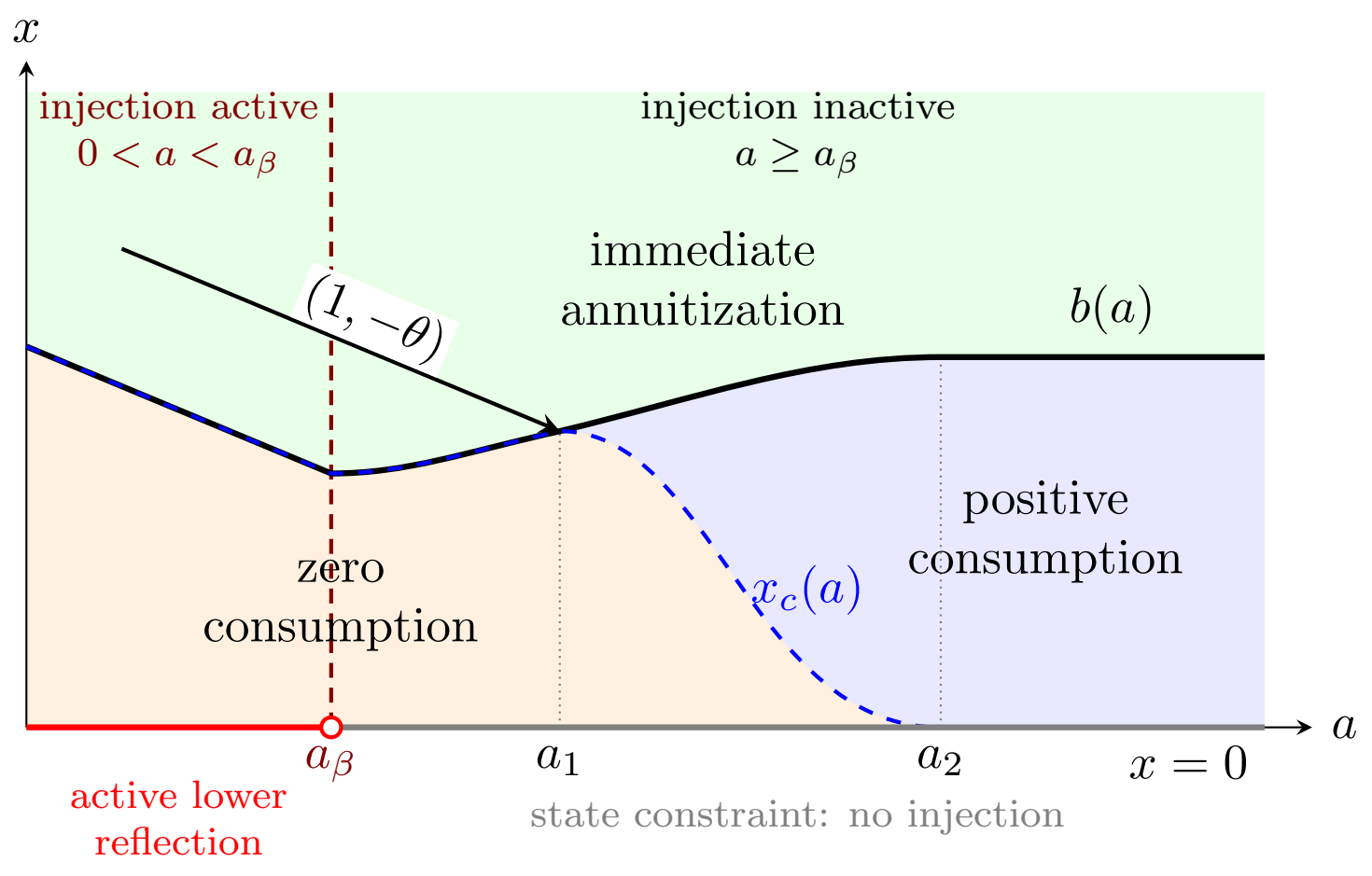
\begin{figure}[H]
\centering
\begin{tikzpicture}[x=1.35cm,y=0.75cm,>=stealth,font=\small]
  \fill[green!9] (0,0) rectangle (6.5,6.0);
  \fill[blue!9] (0,0)--(0,3.6)--(1.6,2.4)
      .. controls (2.0,2.4) and (2.3,2.6) .. (2.8,2.8)
      .. controls (3.5,3.1) and (4.1,3.5) .. (4.8,3.5)
      --(6.5,3.5)
      --(6.5,0)--cycle;
  \fill[orange!12] (0,0)--(0,3.6)--(1.6,2.4)
      .. controls (2.0,2.4) and (2.3,2.6) .. (2.8,2.8)
      .. controls (3.6,2.8) and (3.8,0) .. (4.8,0)--cycle;

  \draw[->] (0,0)--(6.75,0) node[right] {$a$};
  \draw[->] (0,0)--(0,6.3) node[above] {$x$};
  \draw[gray,densely dotted] (2.8,0)--(2.8,2.8);
  \draw[gray,densely dotted] (4.8,0)--(4.8,3.5);
  \draw[thick,densely dashed,DRed] (1.6,0)--(1.6,6.0);

  \draw[very thick] (0,3.6)--(1.6,2.4)
      .. controls (2.0,2.4) and (2.3,2.6) .. (2.8,2.8)
      .. controls (3.5,3.1) and (4.1,3.5) .. (4.8,3.5)
      --(6.5,3.5);
  \draw[thick,blue,dashed] (0,3.6)--(1.6,2.4)
      .. controls (2.0,2.4) and (2.3,2.6) .. (2.8,2.8)
      .. controls (3.6,2.8) and (3.8,0) .. (4.8,0);
  \draw[very thick,red] (0,0)--(1.6,0);
  \draw[very thick,gray] (1.6,0)--(6.5,0);
  \draw[red,fill=white,thick] (1.6,0) circle (2pt);

  \node[DRed,align=center,font=\scriptsize] at (0.8,5.6)
      {injection active\\$0<a<a_\beta$};
  \node[align=center,font=\scriptsize] at (4.05,5.6)
      {injection inactive\\$a\geq a_\beta$};
  \node[DRed,below] at (1.6,0) {$a_\beta$};
  \node[below] at (2.8,0) {$a_1$};
  \node[below] at (4.8,0) {$a_2$};
  \node[below] at (6.1,0) {$x=0$};
  \node[red,align=center,font=\scriptsize] at (0.8,-0.85)
      {active lower\\reflection};
  \node[gray,align=center,font=\scriptsize] at (4.05,-0.85)
      {state constraint: no injection};

  \node[blue] at (4.1,1.3) {$x_c(a)$};
  \node at (5.7,3.95) {$b(a)$};
  \node[align=center] at (1.65,1.15) {zero\\consumption};
  \node[align=center] at (5.35,1.85) {positive\\consumption};
  \node[align=center] at (3.55,4.25) {immediate\\annuitization};
  \draw[->,thick] (0.5,4.525)--(2.8,2.8)
      node[pos=0.55,above,sloped,fill=white,inner sep=1pt]
      {$(1,-\theta)$};
\end{tikzpicture}
\caption{Policy regions and boundaries in the $(a,x)$-plane for
$a_\beta<a_1<a_2$ ($\beta>z_{0b}$).
}
\label{fig:policy-regions}
\end{figure}

\endgroup

\subsubsection{Extension to the annuitization region}

We have constructed the candidate value function on the continuation region $\mathcal C$ through the dual solution. We now extend it to the annuitization region $\mathcal A$. To this end, we first establish a key property of the annuitization boundary.
\begin{lemma}\label{lem:annuity-structure}
Consider annuitization boundary $b$ given by \eqref{eq:positive-boundary-recovery-ba} and define $G(a):=b(a)+\theta a$ for $a>0$. Then  $ G'(a)=0$ for $0<a<a_\beta$, whereas
$ G'(a)>0$ for $a>a_\beta$. Moreover, $ G(a)\to\infty$ as $a\to\infty$.  
\end{lemma}

We next use Lemma~\ref{lem:annuity-structure} to determine the immediate annuity purchase. Let $x>b(a)$ and define $w:=x+\theta a$. Since $x>b(a)$, $w=x+\theta a>b(a)+\theta a=G(a)$. We seek an annuity level $\bar a:=a+\Delta a>a$ such that
\begin{align}\label{eq:jump-annuity}
x-\theta\Delta a=b(a+\Delta a)=b(\bar a),
\end{align}
or equivalently, $G(\bar a)=w$. By Lemma \ref{lem:annuity-structure}, $G$ is strictly increasing on $[a_\beta,\infty)$ and
$G(a)\to\infty$ as $a\to\infty$. Hence, there exists a unique $\bar a>\max\{a,a_\beta\}$ such that $G(\bar a)=w$ and then we set $\bar{x}=b(\bar a)$. The purchase therefore moves the state to the part of the
annuitization boundary where capital injection is inactive.


 We extend the value function to the annuitization region by
\begin{align}\label{eq:annuitization-value-extension}
v(x,a)=v(\bar x,\bar a),\qquad x>b(a).
\end{align}
This value-matching condition corresponds to an immediate
purchase of $\Delta a=\bar a-a$ units of annuity income
at a total cost of $\theta\Delta a$.

For $a<a_\beta$, \eqref{eq:active-annuity-boundary} gives
$b'(a)=-\theta$ along the annuitization boundary. The smooth-fit condition
\eqref{eq:primal-smooth-fit} therefore yields
\[
\frac{d}{da}v(b(a),a)
=y_b(a)\bigl(b'(a)+\theta\bigr)=0.
\]
Thus, the value function is constant along this boundary segment.
An immediate annuity purchase that moves the state to the
junction point $(b(a_\beta),a_\beta)$ therefore leaves the value
unchanged. Section~\ref{sec:verification} specifies this jump
and the subsequent local upper reflection along the part
of the boundary where capital injection is inactive.

\begin{theorem}\label{thm:positive-rate-candidate}
The dual construction and the extension
\eqref{eq:annuitization-value-extension} yield a bounded
candidate value function $v(x,a)$ that is increasing in both state
variables $(x,a)$ and concave in wealth $x$, with strict concavity
on each smooth region.
Furthermore, the candidate value function $v(x,a)$ satisfies HJB-VI \eqref{eq:HJB-v}.
The candidate value function and its first derivatives match across the
intervention boundaries. The second derivatives at these
boundaries are understood as the appropriate one-sided
limits from the adjacent smooth regions.
\end{theorem}

With $L^X\equiv0$ and $\rho=r$, the state equation, discount rate,
consumption constraint, and boundary conditions coincide with those
of \citet{LiangYoung18}. For any fixed $a>0$, this benchmark is also
recovered from the problem with capital injection as $\beta\to\infty$.
Indeed, \eqref{eq:abeta-positive} implies $a_\beta\to0$, so the
upper bound on the marginal value of wealth eventually ceases to bind.

\begin{corollary}[No-injection reduction to \citet{LiangYoung18}]\label{prop:LY-reduction}
If capital injection is prohibited and $\rho=r$, the no-injection
construction above recovers the solution in
\citet[Section~3]{LiangYoung18}.
In particular, the three consumption ranges, the two policy boundaries,
and the regular feedback rules coincide. For any fixed $a>0$, the
solution with capital injection also coincides with this no-injection
solution when $\beta$ is sufficiently large.
\end{corollary}


\section{Verification and Optimal Policies}\label{sec:verification}
In this section, we verify that the candidate constructed from the free-boundary problem coincides with the value function and the associated feedback strategy is optimal. We first state the verification theorem.

\begin{theorem}[Verification Theorem]\label{thm:verification}
Let $v(x,a)$ be the candidate value function defined by
\eqref{eq:inverse-transform} on the continuation region and extended to the
annuitization region by \eqref{eq:annuitization-value-extension}. Then the candidate value function $v(x,a)$ equals $V(x,a)$ given by \eqref{eq:u1}, that is,  $v(x,a)=V(x,a)$, for
$(x,a)\in\mathbb{R}_+^2$.

Furthermore, let $\pi^*(x,a)$ and $c^*(x,a)$ be the  feedback controls given by
\eqref{eq:dual-feedback-detailed} and
\eqref{eq:investment-branches-explicit}. Let $A^*=(A_t^*)_{t\geq0}$ be the minimal non-decreasing c\`adl\`ag
annuitization process that keeps the controlled state in the region $X_t^*\leq b(A_t^*)$, for $t\geq0$. Its continuous part $A^{*,c}$ increases only on the annuitization boundary, in the sense that
$\int_0^\infty
\idc_{\{X_t^*<b(A_t^*)\}}
\,dA_t^{*,c}=0$, and its jumps are determined by \eqref{eq:jump-annuity}. Define the regular
control processes by $\pi_t^*=\pi^*(X_{t-}^*,A_{t-}^*)$,
$c_t^*=c^*(X_{t-}^*,A_{t-}^*)$, where $X^*$ denotes the corresponding controlled wealth process. Then the feedback strategy
$(\pi^*,c^*,A^*)=((\pi_t^*)_{t\geq0},(c_t^*)_{t\geq0},(A_t^*)_{t\geq0})\in\mathbb U$ is optimal.
\end{theorem}

\begin{figure}[H]
\centering
\begin{tikzpicture}[x=1.25cm,y=0.9cm,>=stealth,font=\small]

\fill[red!4]  (0.45,0.45) rectangle (3.30,4.85);
\fill[blue!3] (3.30,0.45) rectangle (8.35,4.85);

\draw[->] (0.45,0.45)--(8.55,0.45) node[right] {$a$};
\draw[->] (0.45,0.45)--(0.45,5.10) node[above] {$x$};

\draw[densely dashed,gray] (3.30,0.45)--(3.30,4.85);
\node[below] at (3.30,0.45) {$a_\beta$};

\draw[gray,densely dotted] (5.10,0.45)--(5.10,2.95);
\node[below] at (5.10,0.45) {$a_1$};

\draw[gray,densely dotted] (6.70,0.45)--(6.70,3.65);
\node[below] at (6.70,0.45) {$a_2$};

\node[align=center] at (1.90,5.2)
    {injection active\\[-1pt]$0<a<a_\beta$};

\node[align=center] at (5.95,5.2)
    {injection inactive\\[-1pt]$a\ge a_\beta$};

\draw[very thick,blue]
    (0.80,3.60)--(3.30,2.40)
    .. controls (4.00,2.40) and (4.45,2.60) .. (5.10,2.80)
    .. controls (5.85,3.10) and (6.25,3.50) .. (6.70,3.50)
    -- (8.00,3.50);

\node[font=\scriptsize,blue,anchor=west] at (7.85,3.68)
    {$x=b(a)$};

\coordinate (AB) at (3.30,2.40);   


\coordinate (D1) at (2.20,4.18);
\coordinate (J1) at (5.10,2.80);

\coordinate (D2) at (1.35,1.70);
\coordinate (B2) at (1.35,3.34);   
\coordinate (R2) at (4.55,2.64);   

\coordinate (D3) at (6.10,4.15);
\coordinate (B3) at (7.00,3.50);   

\coordinate (D4) at (5.55,1.90);
\coordinate (B4) at (5.55,2.99);   
\coordinate (R4) at (6.35,3.38);   

\fill (D1) circle (1.9pt);
\fill (D2) circle (1.9pt);
\fill (D3) circle (1.9pt);
\fill (D4) circle (1.9pt);

\node[font=\scriptsize,above left]  at (D1) {$D_1$};
\node[font=\scriptsize,left]        at (D2) {$D_2$};
\node[font=\scriptsize,above left]  at (D3) {$D_3$};
\node[font=\scriptsize,below left]  at (D4) {$D_4$};

\fill (J1) circle (1.5pt);
\fill (B2) circle (1.2pt);
\fill (AB) circle (1.5pt);
\fill (B3) circle (1.5pt);
\fill (B4) circle (1.2pt);

\draw[->,very thick,orange!85!black]
    (D1)--(J1);

\draw[->,densely dashed,thick,gray!80]
    (D2)--(B2);

\draw[->,very thick,orange!85!black]
    (B2)--(AB);

\draw[->,very thick,orange!85!black]
    (AB) .. controls (3.85,2.42) and (4.20,2.52) .. (R2);

\draw[->,very thick,orange!85!black]
    (D3)--(B3);

\draw[->,densely dashed,thick,gray!80]
    (D4)--(B4);

\draw[->,very thick,orange!85!black]
    (B4) .. controls (5.80,3.10) and (6.05,3.28) .. (R4);


\node[font=\scriptsize,orange!85!black,rotate=-15] at (3.95,3.55)
    {annuity purchase};
\end{tikzpicture}

\caption{Illustration of the annuity purchase process.}
\label{fig:annuity-purchase-geometry}
\end{figure}
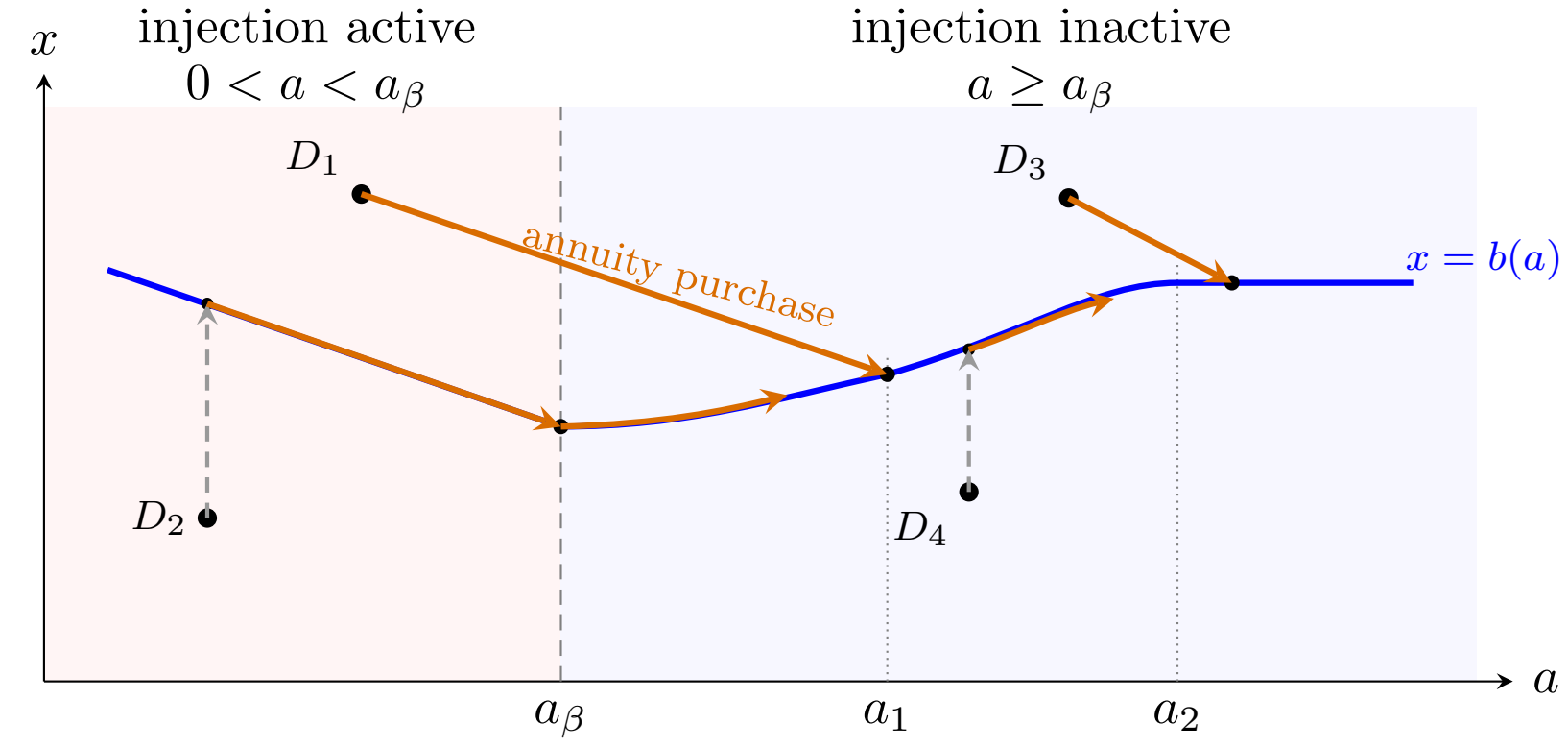

 By Theorem \ref{thm:verification}, the optimal annuitization process is the minimal non-decreasing intervention process that keeps the state in the region $X_t^*\leq b(A_t^*)$ and its continuous part satisfies $\int_0^\infty
    \idc_{\{X_t^*<b(A_t^*)\}}
    dA_t^{*,c}=0$, while its jumps are determined by the projection rule \eqref{eq:jump-annuity}. Unlike the classical annuitization models without capital injection
in \cite{MilevskyYoung07,LiangYoung18}, where the optimal strategy is typically characterized by a pure reflection structure, the presence of
capital injection leads to a mixed jump--reflection structure. The qualitative behavior depends on whether the state lies on the injection-active or injection-inactive part of the annuitization
boundary, as illustrated in Figure~\ref{fig:annuity-purchase-geometry}.

On the injection-inactive part of the boundary, where
$b'(a)>-\theta$, annuity purchases are implemented through the standard
continuous Skorokhod reflection. Specifically, when the state is above
the boundary (as illustrated by $D_3$ in Figure \ref{fig:annuity-purchase-geometry}), the optimal policy immediately
purchases annuity income to project the state onto $x=b(a)$. Meanwhile, for a
state below the boundary (as illustrated by $D_4$), wealth accumulation
first drives the state to the boundary, after which annuity purchases
continue along the boundary through continuous reflection.  If it reaches $x=b(a)$, annuity purchases are then made as needed through continuous reflection to keep the state in the region $x\le b(a)$.

The behavior is fundamentally different on the injection-active part of the boundary. In this region, the boundary satisfies $b'(a)=-\theta$, and therefore the annuity-purchase direction is tangent to the free boundary. Consequently,
a marginal annuity purchase only transfers liquid wealth into irreversible annuity income and the continuous annuitization cannot keep the state in the region $x\le b(a)$. For a state above the boundary (as illustrated by $D_1$), the optimal
strategy therefore performs an immediate discrete annuity purchase that
moves the state across the boundary point $a_{\beta}$ into the
injection-inactive part of the boundary. For a state initially below the boundary (as illustrated by $D_2$), the state evolves under the regular controls without annuity purchases and may reach the injection-active part of the boundary. If it does, an immediate annuity purchase can move the state along the boundary to the point $(b(a_\beta),a_\beta)$. Once the state enters the injection-inactive region, further annuitization is again
implemented through continuous reflection along the annuitization
boundary to keep the state in the region $x\le b(a)$.

%
%

\section{Properties of the Optimal Strategies and Numerical Examples}\label{sec:numerical}

In this section, we present some numerical examples to illustrate our theoretical fings. We compare the optimal strategies with their benchmark in
\citet{LiangYoung18}  when capital injections are not unavailable. Furthermore, we investigate the sensitivity of the optimal policies with
respect to key model parameters and simulate the corresponding wealth and
intervention processes to highlight the financial implications of the proposed strategies.

Throughout the numerical experiments, we adopt the parameter values in
Table~\ref{tab:calibration-positive} as the benchmark calibration. The
parameters are chosen to be close to those used in \citet{LiangYoung18},
while incorporating the additional capital injection feature considered in
our model.
\begin{table}[H]
\centering
\small
\begin{tabular}{ccl}
\hline
Parameter & Value & Interpretation\\
\hline
$r$ & $0.04$ & risk-free interest rate\\
$\rho$ & $0.04$ & pure time-discount rate\\
$\lambda^s$ & $0.06$ & subjective mortality intensity\\
$\lambda^o$ & $0.04$ & objective mortality intensity\\
$\mu$ & $0.08$ & risky asset expected return\\
$\sigma$ & $0.20$ & risky volatility\\
$\gamma$ & $0.02$ & absolute risk aversion \\
$\beta$ & $1.20$ & marginal cost of capital injection\\
\hline
\end{tabular}
\caption{Model parameters used in the benchmark numerical example.}
\label{tab:calibration-positive}
\end{table}


\subsection{Pathwise illustration of the optimal strategy}\label{subsec:pathwise-simulation}
To visualize the interaction among investment, consumption, annuitization, and capital injection, we examine a representative sample path generated under the benchmark calibration in Table~\ref{tab:calibration-positive}.  We simulate the controlled state
process under the optimal investment, consumption, and annuitization
strategies by Monte Carlo method, with initial state $X_0^*=5, A_0^*=1$. For the pathwise analysis, define the stopping time
$\tau_\beta:=\inf\{t\geq0:A_t^*\geq a_\beta\}$, namely, the first time at which annuity income reaches the threshold
$a_\beta$. Before $\tau_\beta$, the state lies in the injection-active
regime, whereas after $\tau_\beta$ the lower wealth boundary becomes a state constraint and no further capital injection is required.

\begin{figure}[H]
\centering

 \subfigure[]{
\includegraphics[width=0.45\textwidth]{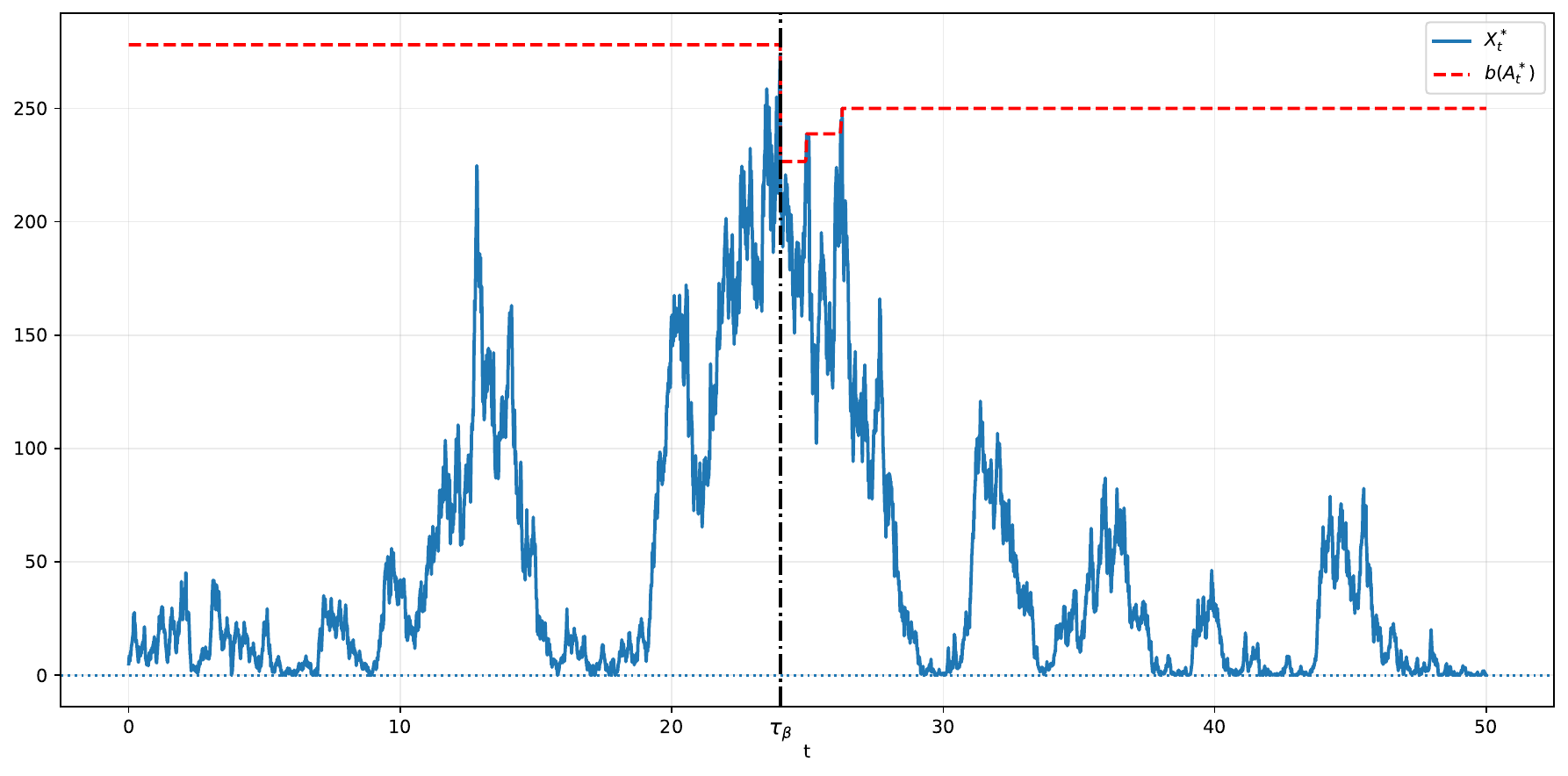}
}
\subfigure[]{
\includegraphics[width=0.45\textwidth]{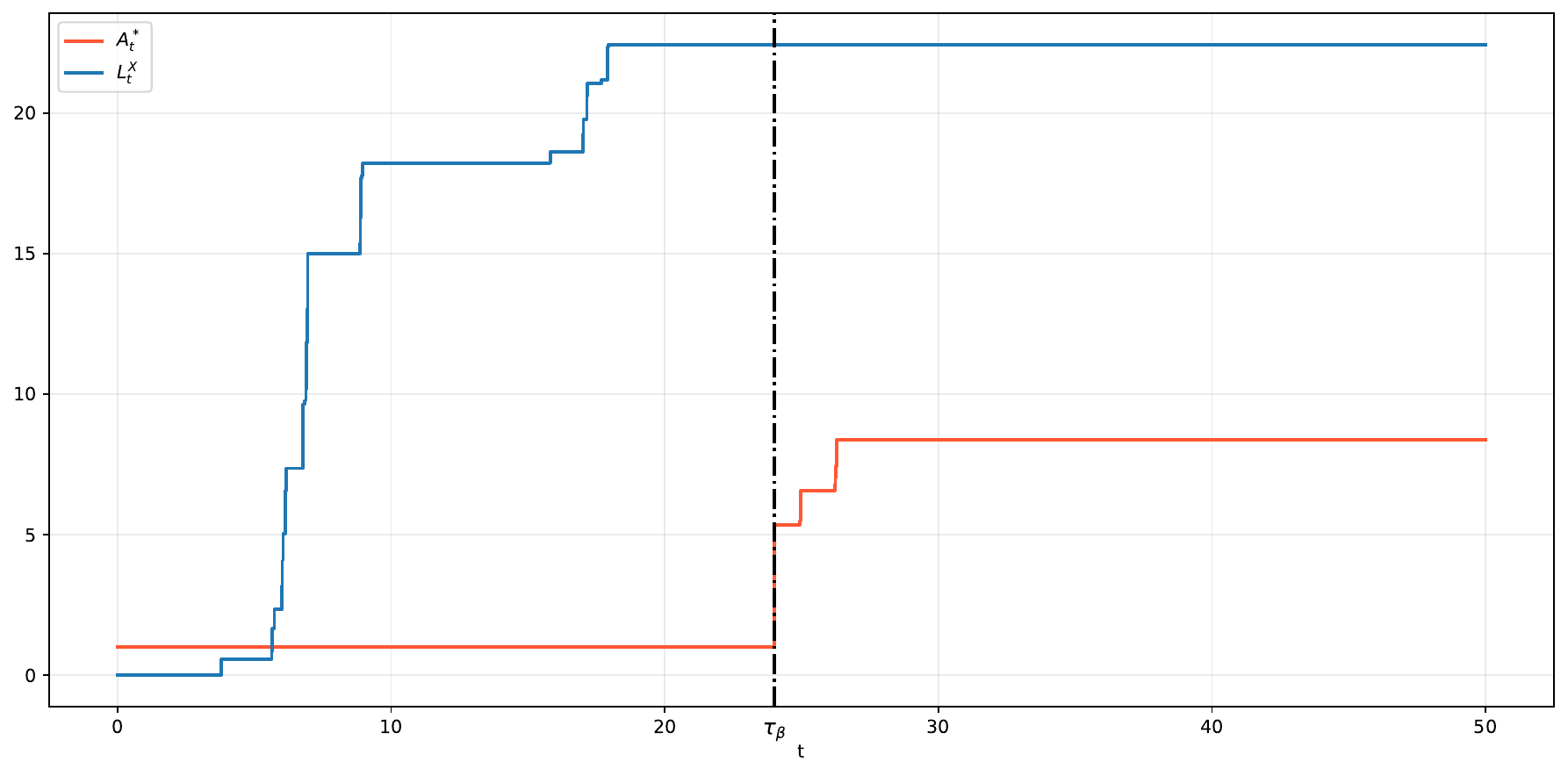}
}
\caption{(a): A sample path of optimal liquid wealth $X_t^*$ together with the annuitization boundary $b(A_t^*)$. (b): Optimal annuity income $A_t^*$ and  capital injection
$L_t^{X}$ along the same sample path.  }
\label{fig:pathwise-wealth}
\end{figure}

Figure~\ref{fig:pathwise-wealth}-(a) plots the optimal liquid wealth process
$X_t^*$ together with the corresponding annuitization boundary
$b(A_t^*)$, while Figure~\ref{fig:pathwise-wealth}-(b) displays the
optimal annuity process $A_t^*$ and the capital injection
process $L_t^{X,*}$. While $A_t^*<a_\beta$, capital injection occurs only when liquid wealth reaches the lower boundary $X_t^*=0$. Along the displayed path, several injections occur each time wealth reaches zero, illustrating the role of the solvency guarantee in supporting the account against adverse market realizations. Throughout this period, no annuity is purchased because liquid wealth remains below the annuitization boundary.

A qualitatively different intervention occurs at $\tau_\beta$. At this
time, the annuity process makes a discrete upward adjustment from its
pre-jump level to a level above $a_\beta$, thereby moving the state
directly to the injection-inactive part of the annuitization boundary,
while liquid wealth is reduced by the corresponding annuity premium.
Subsequent interactions with the upper boundary generate several smaller annuity purchases. In contrast, the  capital injection process $L^X_t$ remains constant after $\tau_\beta$. The simulated path therefore illustrates the change in the optimal policy at the endogenous switching level $a_\beta$. Before $\tau_\beta$, maintaining solvency may require capital injection when liquid wealth reaches zero, whereas after $\tau_\beta$ the account remains solvent without further external support and annuity purchases are determined by the upper boundary.

\subsection{Comparison of optimal strategy with the no-injection benchmark}
\label{subsec:comparison-no-injection}

We next examine how the presence of capital injections affects the
optimal annuitization, consumption, and investment policies. Under the setup in Table~\ref{tab:calibration-positive}, we compare the optimal policies in our model with the no-injection benchmark of \citet{LiangYoung18}.  We use the superscripts I and N to denote the policies with and without capital injection, respectively.

\begin{figure}[ht]
\centering
 \subfigure[]{
        \includegraphics[width=0.3\textwidth]{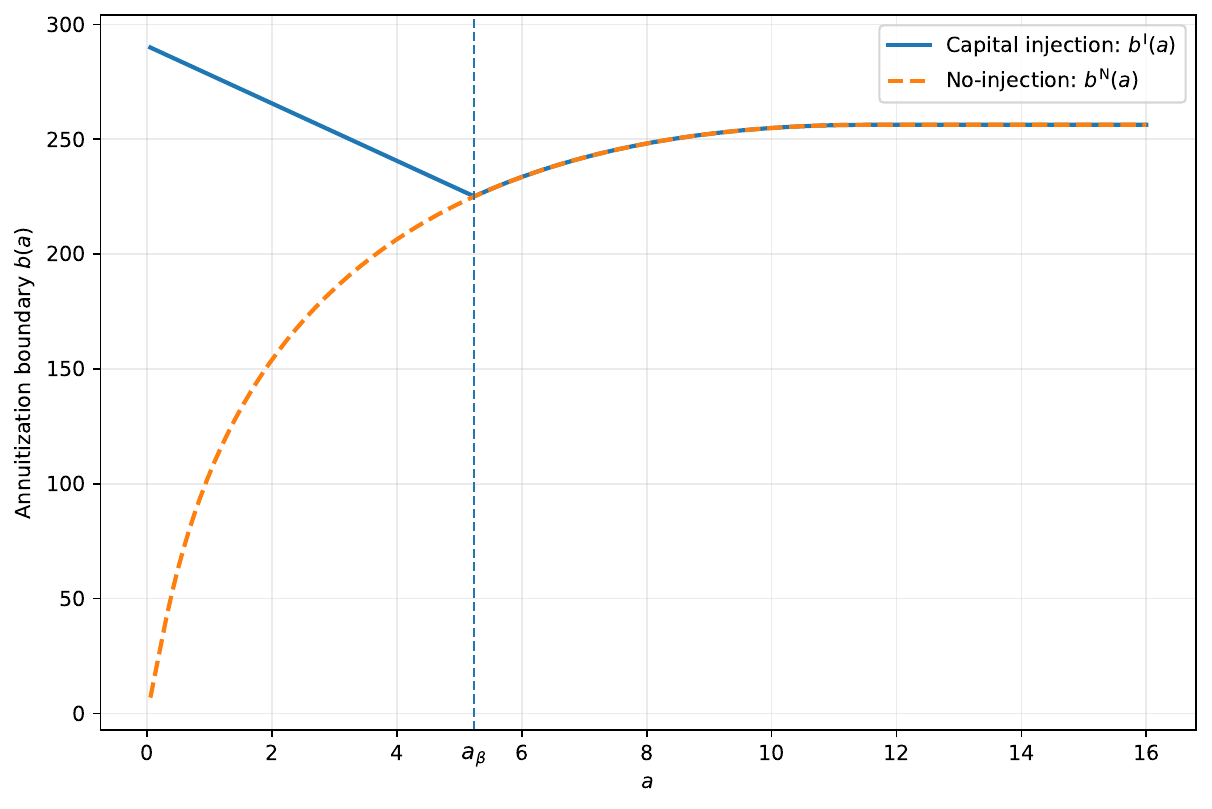}
    }
     \subfigure[]{
        \includegraphics[width=0.3\textwidth]{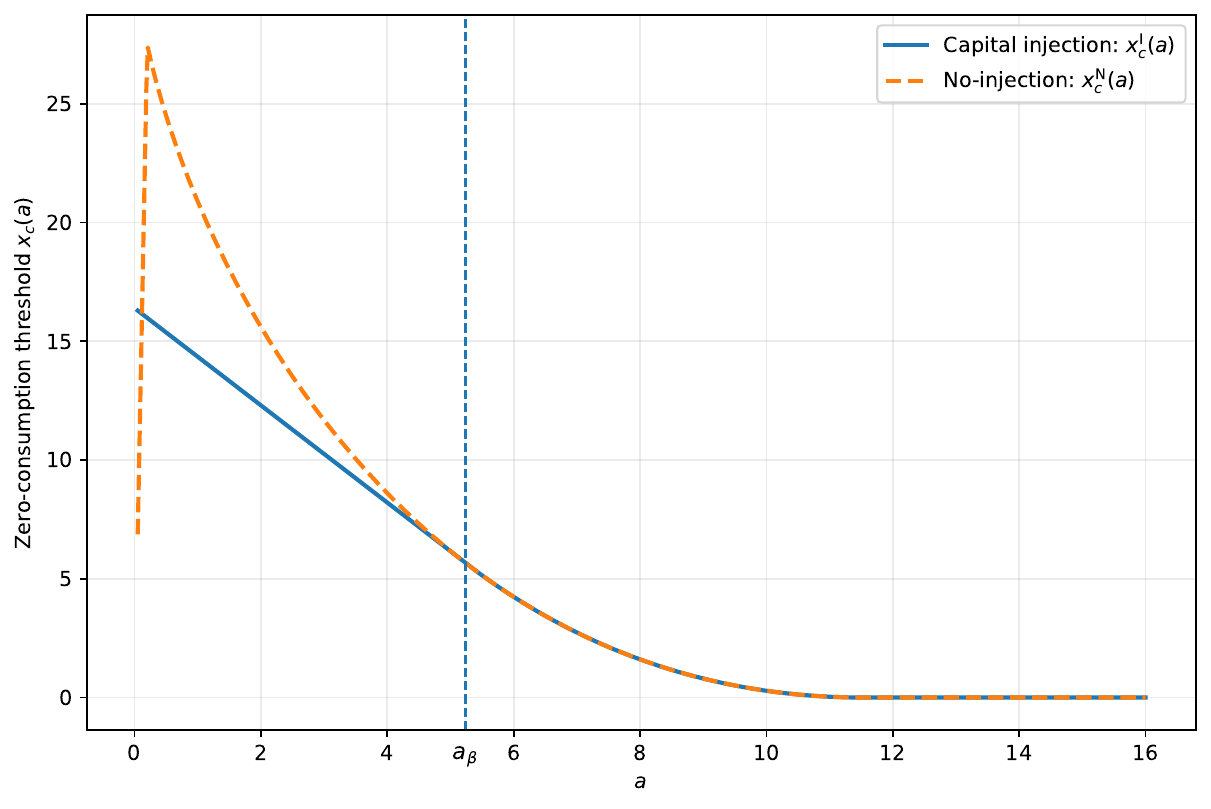}
    }
         \subfigure[]{
        \includegraphics[width=0.3\textwidth]{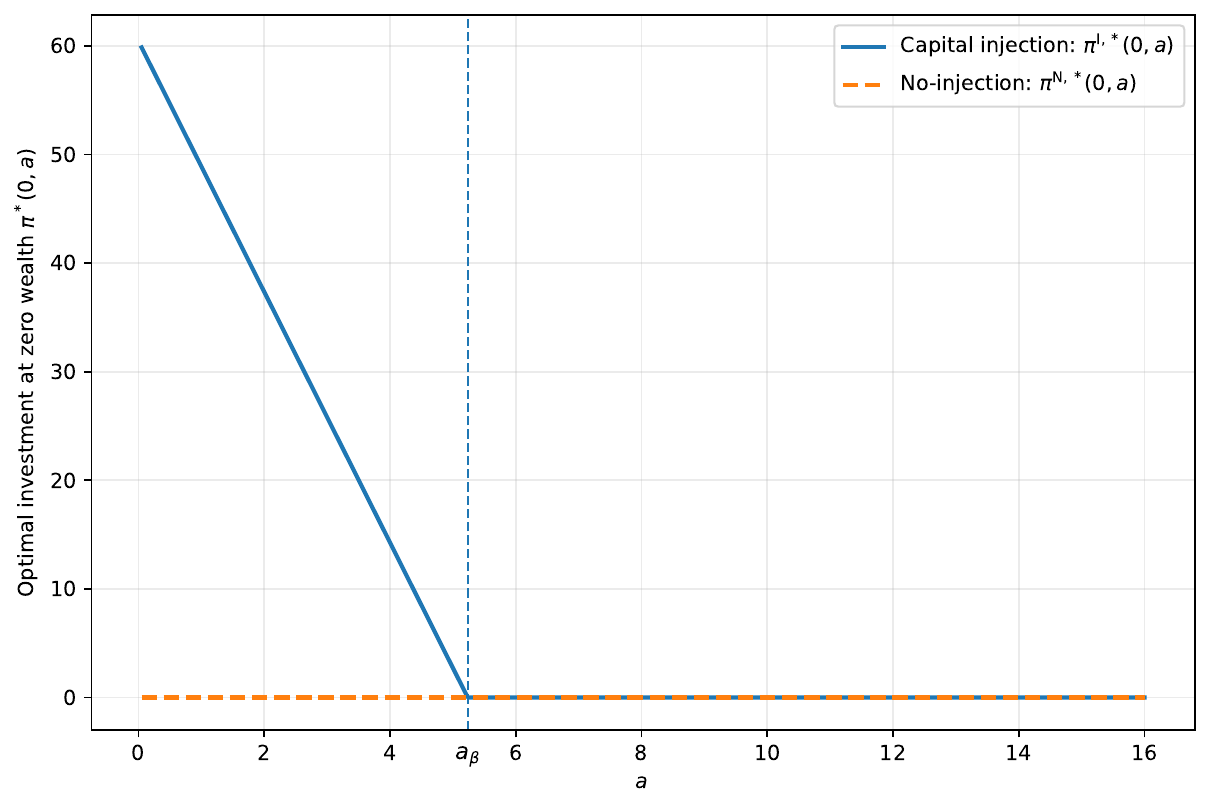}
    }
\caption{Comparison of the optimal policies with and without capital
injection under the benchmark calibration. Panels (a)--(c) show the
annuitization boundary $b(a)$, the zero-consumption boundary $x_c(a)$, and the optimal risky investment $\pi^*(0,a)$ at zero liquid wealth,
respectively. }
\label{fig:policy-comparison}
\end{figure}

Figure~\ref{fig:policy-comparison} compares the annuitization boundary, the zero-consumption boundary, and the optimal risky investment at zero liquid wealth. The vertical line at $a_\beta$ separates the injection-active and injection-inactive regions. As $\rho=r$ in the model setup, Corollary~\ref{prop:LY-reduction} implies that the no-injection solution derived in Section~\ref{sec:dual} coincides with that of \citet{LiangYoung18}. As shown in Figure~\ref{fig:policy-comparison}, the two models yield identical optimal policies for $a\geq a_\beta$, so the effects of capital injection are confined to the region $a<a_\beta$.

Figure~\ref{fig:policy-comparison}-(a) reveals two noteworthy and somewhat counterintuitive features of the optimal annuitization policy. First, access to capital injection raises the annuitization boundary when
$a<a_\beta$, that is, $b^{\text{I}}(a)>b^{\text{N}}(a)$, for  $a<a_\beta$. Hence, when solvency support is available, the pension-plan manager optimally retains a larger amount of liquid wealth before purchasing additional annuity income. At first sight, one might expect the availability of external support to encourage earlier annuitization. The opposite result arises because an annuity purchase is irreversible and converts $\theta$ units of liquid wealth into one unit of future annuity income. Capital injection provides downside protection at the lower wealth boundary and therefore reduces the need to use annuitization itself as a device for protecting future consumption against low-wealth states. At the same time, retaining liquid wealth preserves portfolio flexibility and allows the manager to continue taking advantage of risky investment opportunities. Consequently, the solvency guarantee increases the option value of remaining in the continuation region and delays annuity purchases. 

Second, and more strikingly, the annuitization boundary is decreasing in existing annuity income throughout the injection-active region. 
This behavior is in sharp contrast with the no-injection model of
\citet{LiangYoung18}, where the annuitization boundary is increasing in
existing annuity income before eventually becoming constant. The difference is caused by the solvency guarantee. In the injection-active region, it follows from Lemma \ref{lem:annuity-structure} that
$b^{\text{I}}(a)+\theta a=b^{\text{N}}(a_\beta)+\theta a_\beta$, so the relevant total annuitizable resource at the upper boundary remains constant. An additional unit of existing annuity income therefore reduces the amount of liquid wealth required to reach the annuitization boundary by exactly its purchase price $\theta$. Economically, annuity income and liquid wealth act as substitutes along this part of the boundary. As the retiree accumulates more guaranteed lifetime income, less liquid wealth is required before further annuitization becomes optimal.

These two observations reveal both distinctions and connections between
solvency insurance and longevity insurance. Capital injection protects the retiree against shortfall risk, while annuitization protects future
consumption against longevity risk. When annuity income is low, the
solvency guarantee is valuable and increases the value of retaining
liquidity, thereby pushing the annuitization boundary upward. As annuity income increases toward $a_\beta$, the need for external solvency support gradually disappears, and the liquid-wealth threshold for further annuitization declines. At $a=a_\beta$, the capital-injection constraint ceases to bind, and the funded and no-injection annuitization boundaries join continuously. Thus, the model suggests a form of substitution between sponsor-provided solvency protection and private annuitization that is absent from standard annuitization models such as
\citet{MilevskyYoung07,LiangYoung18}.

Figure \ref{fig:policy-comparison}-(b) compares the wealth threshold
$x_c(a)$ separating the zero- and positive-consumption regions. Capital injection lowers this threshold for $a<a_\beta$, indicating that positive consumption becomes optimal at a lower level of liquid wealth when the solvency guarantee is available. This is because the availability of external support reduces the need to preserve liquid wealth solely as a buffer against future shortfalls, allowing the retiree to begin consuming earlier as wealth increases. Another direct effect of capital injection is illustrated in Figure~\ref{fig:policy-comparison}-(c). In \citet{LiangYoung18}, non-negative wealth is imposed as a hard state constraint, which requires $\pi^{\text{N},*}(0,a)=0$. By contrast, when capital injection is active, the manager can maintain a strictly positive risky position even at zero liquid wealth.

\subsection{Sensitivity analysis}\label{sec:sensitivity-analysis}

\begin{figure}[ht]
\centering

 \subfigure[]{
    \centering
    \includegraphics[width=0.4\textwidth]
    {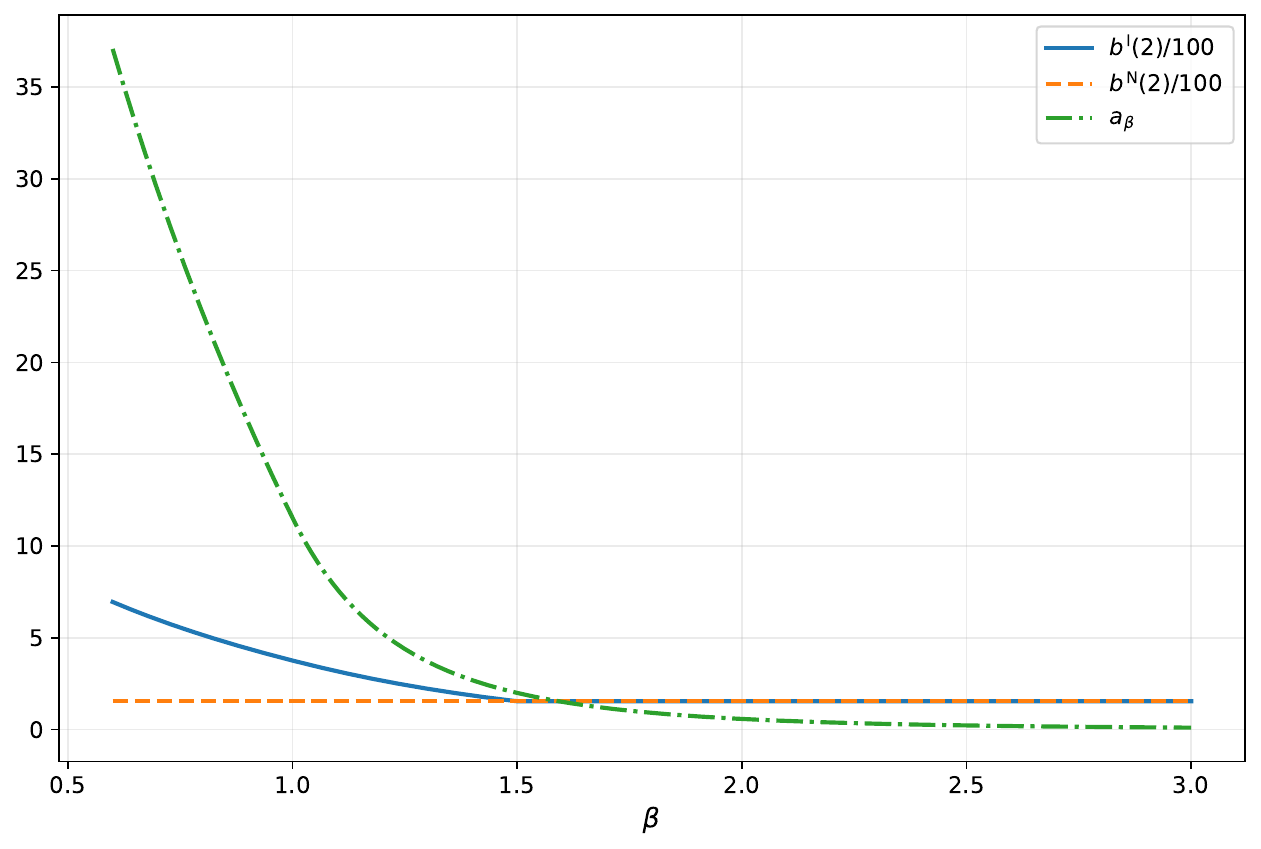}
}
 \subfigure[]{
    \centering
    \includegraphics[width=0.4\textwidth]
    {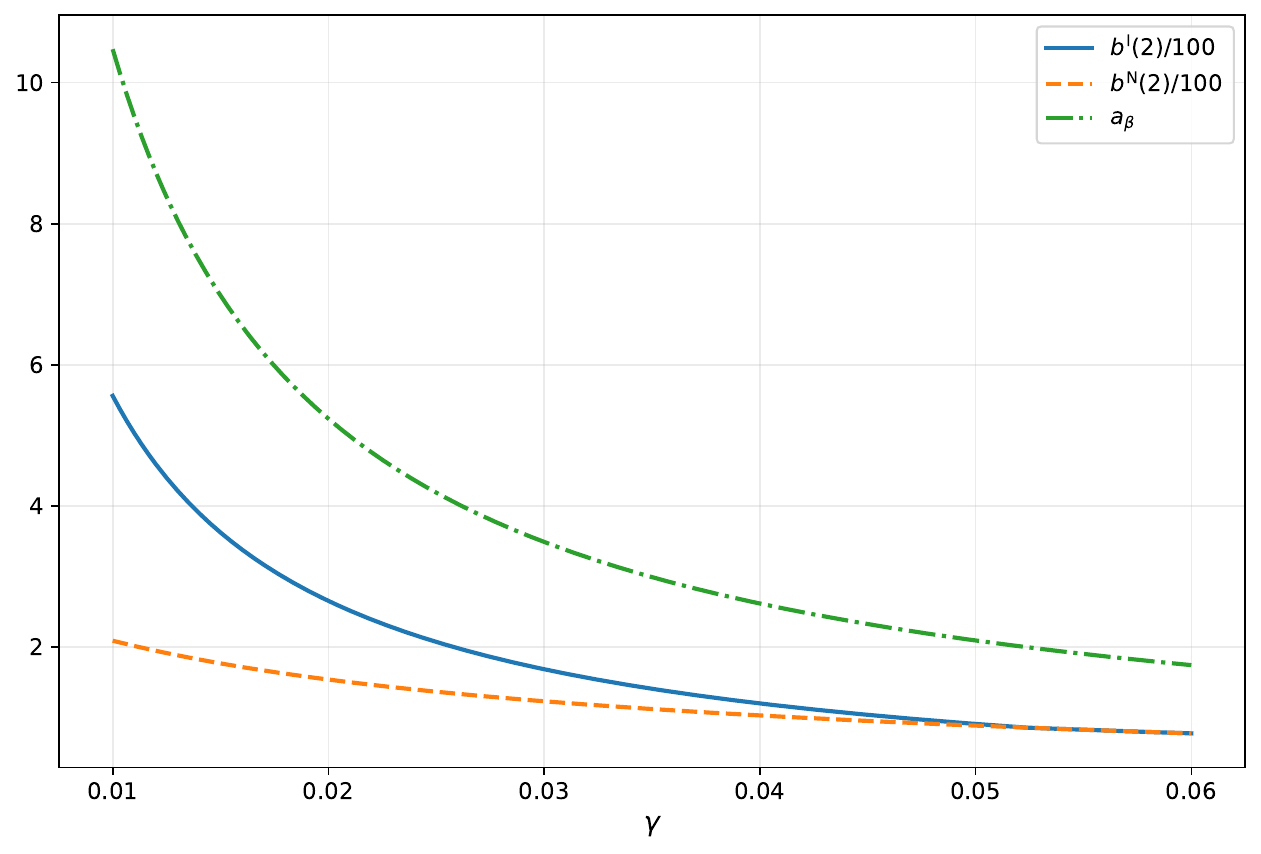}
}

 \subfigure[]{
    \centering
    \includegraphics[width=0.4\textwidth]
    {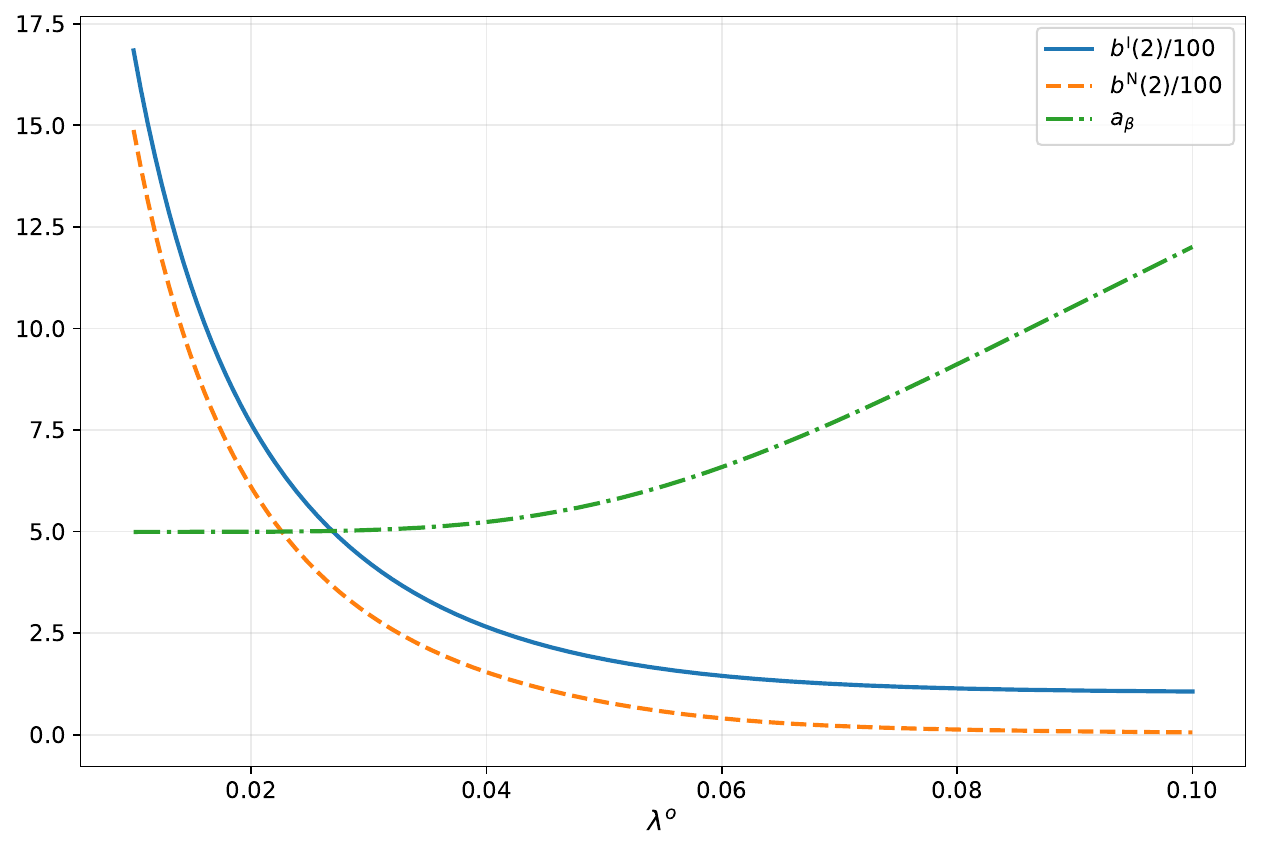}
}
 \subfigure[]{
    \centering
    \includegraphics[width=0.4\textwidth]
    {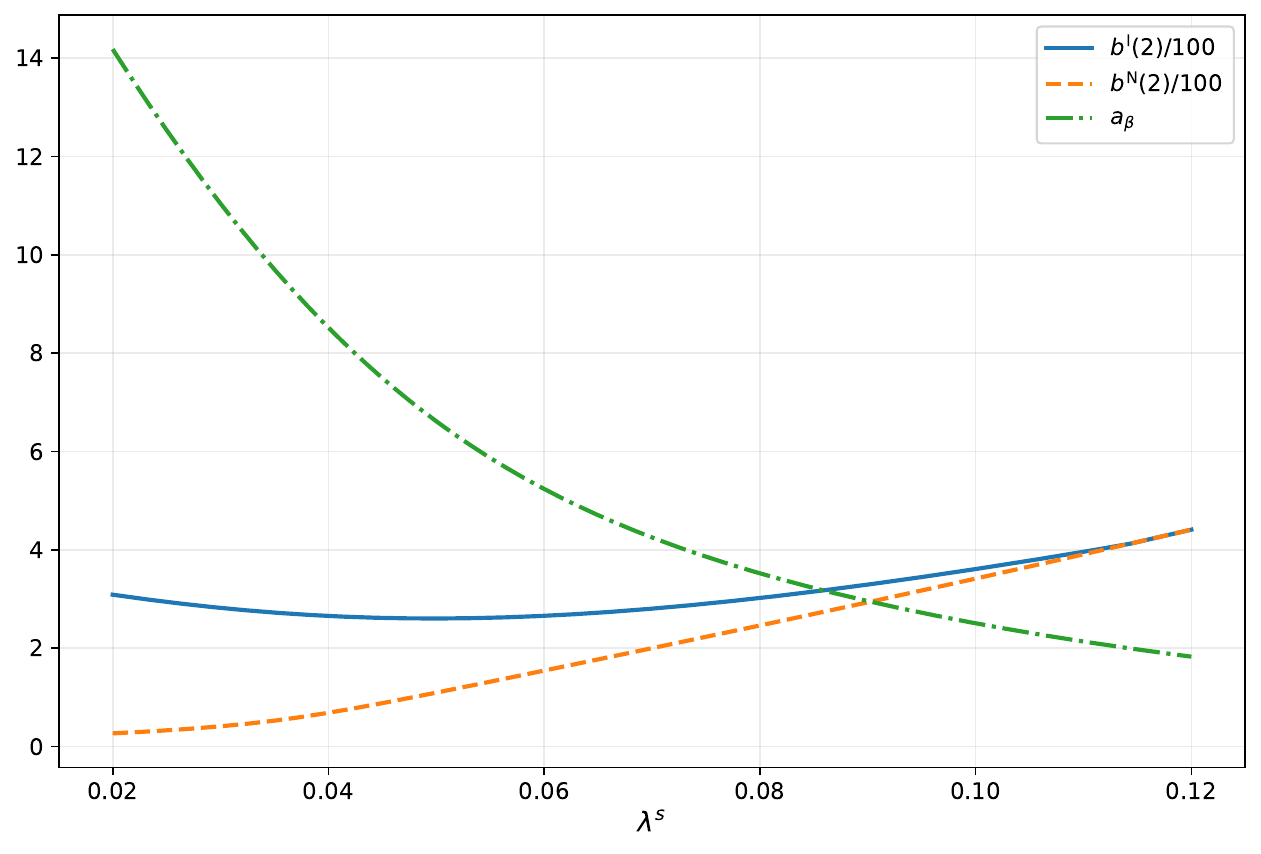}
}

\caption{Sensitivity of the switching annuity level $a_\beta$ and the optimal annuitization boundaries $b^I(a)$ and $b^N(a)$ to key model parameters with $a=2$.}
\label{fig:sensitivity-boundaries}
\end{figure}

\begin{figure}[htbp]
\centering
 \subfigure[]{
    \centering
    \includegraphics[width=0.4\textwidth]
    {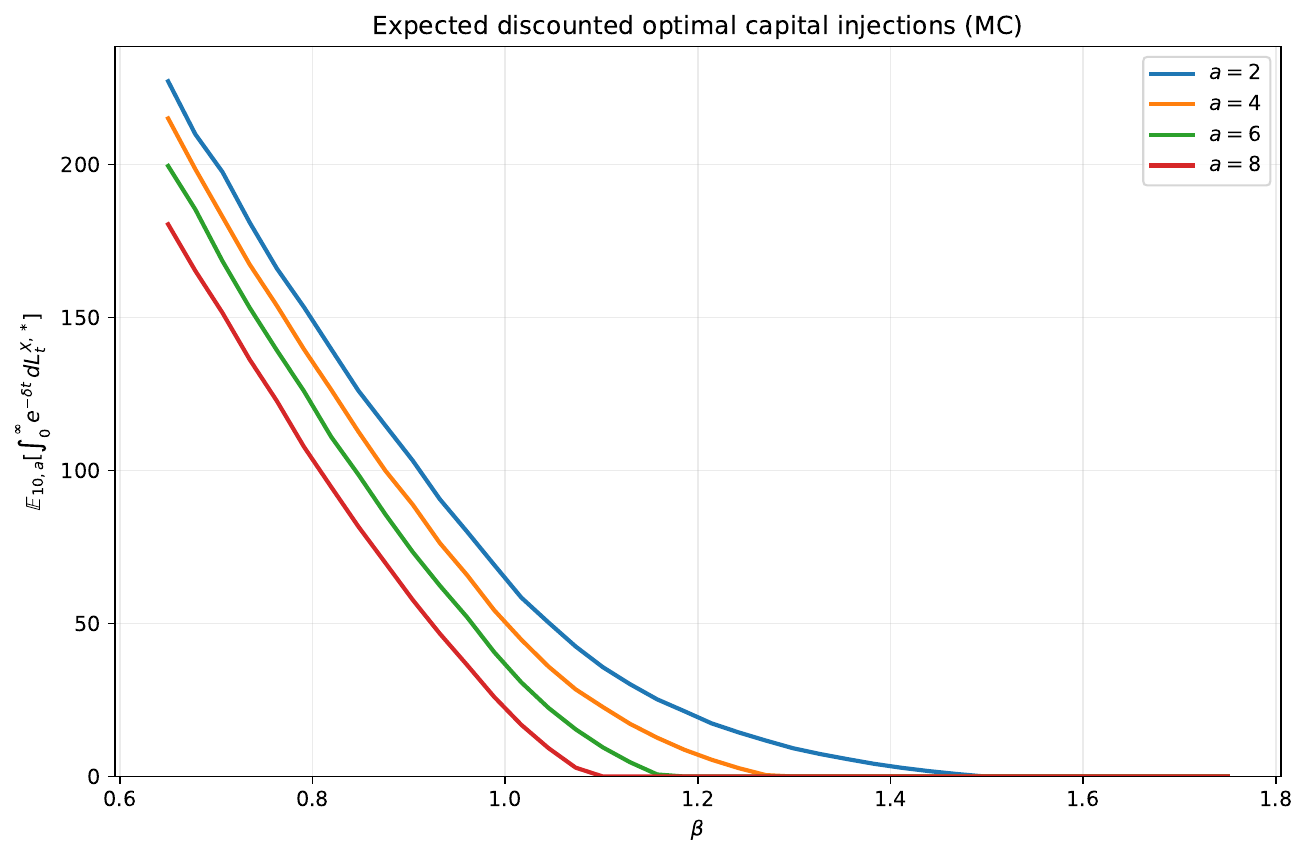}
}
 \subfigure[]{
    \centering
    \includegraphics[width=0.4\textwidth]
    {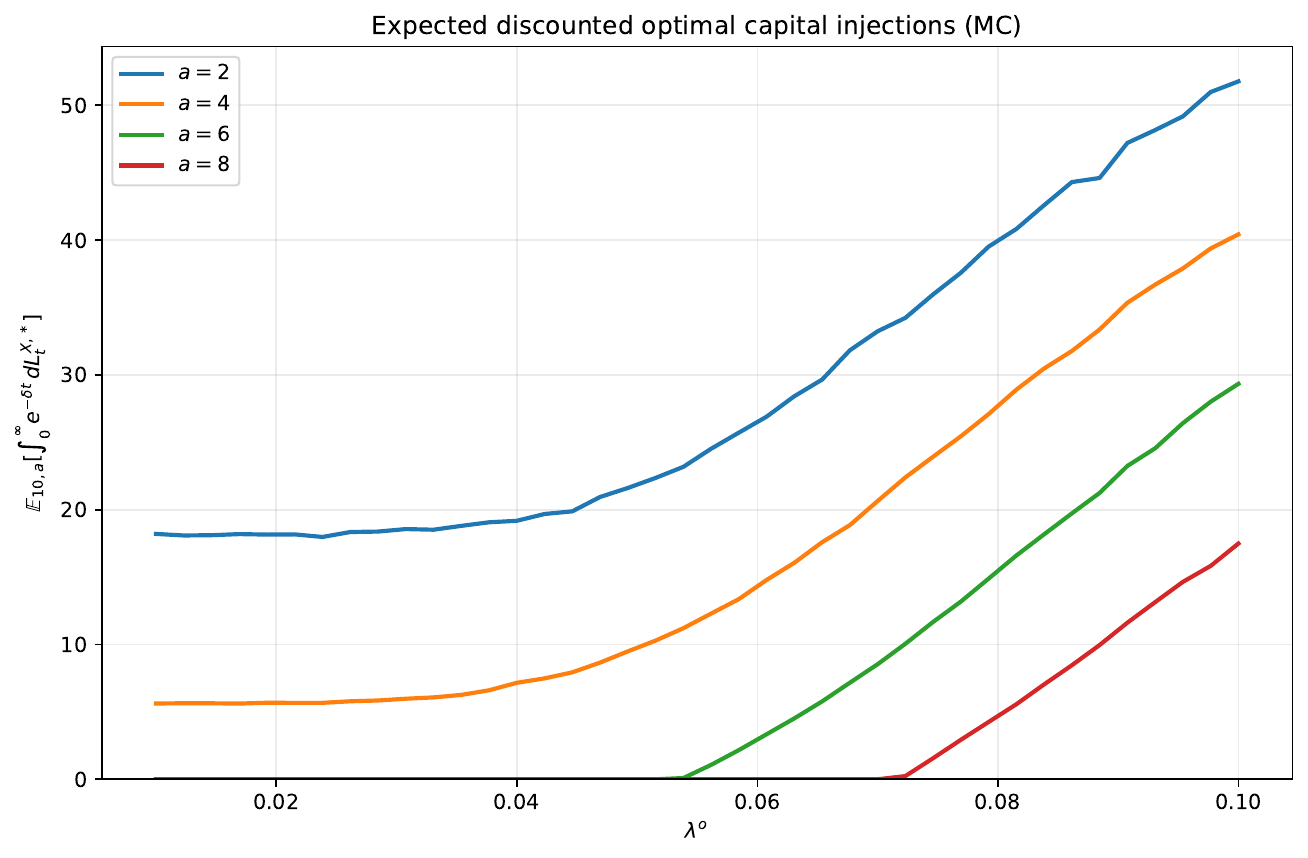}
}
\caption{Sensitivity of the expected discounted cumulative capital injections $\mathbb{E}[\int_0^\infty e^{-\delta t}dL_t^{X}]$
to the marginal injection cost $\beta$ and the objective mortality intensity $\lambda^o$. The initial liquid wealth is $x=10$, and
the initial annuity income takes values $a=2,4,6,8$. }
\label{fig:sensitivity-injections}
\end{figure}

In this subsection, we examine the sensitivity of the optimal annuitization policy and the expected capital injections to key model parameters. We vary one parameter at a time while keeping the remaining parameters fixed as in Table~\ref{tab:calibration-positive}.  Figure~\ref{fig:sensitivity-boundaries} illustrates the sensitivity of the switching level $a_\beta$ and the optimal annuitization boundaries $b^{\text{I}}(a)$ and $b^{\text{N}}(a)$, evaluated at $a=2$. For ease of comparison, the two wealth boundaries are divided by 100 so that all three quantities can be displayed on the same vertical axis. Figure~\ref{fig:sensitivity-injections} plots the expected discounted cumulative capital injection $\mathbb{E}[
\int_0^\infty e^{-\delta t}dL_t^{X}]$ under the optimal policy, with respect to $\beta$ and $\lambda^o$, respectively.

Figure~\ref{fig:sensitivity-boundaries}-(a) shows that both $a_\beta$ and $b^{\text{I}}(a)$ decrease as the marginal cost of capital injection $\beta$ increases. A larger $\beta$ makes external capital support more expensive, reducing the attractiveness of retaining liquid wealth and relying on future injections to maintain solvency. Consequently, the manager optimally purchases additional annuity income at a lower level of liquid wealth. As $\beta$ increases further, the switching level $a_\beta$ eventually falls below the initial annuity income, at which point capital injection becomes inactive and $b^{\text{I}}(a)$ coincides with $b^{\text{N}}(a)$. In contrast, the no-injection boundary is independent of $\beta$, since this parameter is irrelevant when capital injection is prohibited.

Figure~\ref{fig:sensitivity-boundaries}-(b) illustrates the effect of risk aversion. As absolute risk aversion parameter $\gamma$ increases, both annuitization boundaries decrease, indicating that a more risk-averse manager is willing to purchase additional guaranteed annuity income at a lower level of liquid wealth. This is consistent with the findings of \citet{LiangYoung18} in the absence of capital injection. Moreover, the switching level $a_\beta$ also decreases with $\gamma$. Greater risk aversion reduces the willingness to maintain risky investment positions supported by potential future capital injections, thereby narrowing the region in which external support is optimally utilized. For sufficiently large $\gamma$, the two annuitization boundaries coincide at the given annuity income.

Figures~\ref{fig:sensitivity-boundaries}-(c) and (d) distinguish the effects of the objective and subjective mortality intensities. An interesting and somewhat counterintuitive finding emerges from the sensitivity analysis with respect to the objective mortality intensity $\lambda^o$. As $\lambda^o$ increases, the annuity purchase price $\theta=1/(r+\lambda^o)$ decreases, making guaranteed lifetime income more affordable. Consequently, both annuitization boundaries decline, as shown in Figure~\ref{fig:sensitivity-boundaries}-(c), indicating that the manager optimally purchases additional annuity income at a lower level of liquid wealth. One might therefore expect cheaper annuities to reduce the future need for external capital support. However,  the switching level $a_\beta$ increases with $\lambda^o$ and Figure~\ref{fig:sensitivity-injections}-(b) also shows that the expected discounted cumulative capital injections generally increase with $\lambda^o$ under the benchmark calibration. These results reveal an important distinction between the cost of purchasing longevity insurance and the cost of providing solvency protection. A reduction in annuity prices benefits the manager by making guaranteed future income more accessible, but it need not reduce the expected capital contributions required from the pension plan. 

Another interesting finding points to the effect of subjective mortality intensity $\lambda^s$. In the no-injection benchmark, the annuitization boundary increases monotonically with $\lambda^s$, that is,  a shorter expected remaining lifetime reduces the value of future lifetime income. By contrast, Figure~\ref{fig:sensitivity-boundaries}-(d) shows that with capital injection the boundary is non-monotonic. It first decreases and then increases as $\lambda^s$ rises. This reflects two competing effects. On one hand, a higher $\lambda^s$ raises the effective discount rate $\delta=\rho+\lambda^s$, lowering the present value of future annuity payments and encouraging the retiree to retain liquid wealth, which pushes the boundary up. On the other hand, a higher $\lambda^s$ lowers the switching level $a_\beta$, shrinking the range over which external solvency support is optimally used. This reduces the value of retaining liquidity under the guarantee and tends to push the boundary down. Thus, the effect of longevity expectations on annuitization depends on whether external solvency protection is available. In particular, a shorter subjective life expectancy need not delay annuitization when capital injections are part of the retirement financing arrangement; the optimal response balances the lower value of future lifetime income against the diminished benefit of holding liquidity under the solvency guarantee.

Figure~\ref{fig:sensitivity-injections}-(a) shows that the expected discounted capital injections decrease as $\beta$ increases. When external capital is relatively inexpensive, the manager is more willing to maintain a risky investment strategy supported by the pension plan's solvency guarantee, resulting in greater expected future capital contributions. Once the initial annuity income satisfies $a\geq a_\beta$, no further capital injection is required under the optimal strategy, and the corresponding expectation becomes zero. Figure~\ref{fig:sensitivity-injections}-(b) reveals a different pattern with respect to $\lambda^o$.  As $\lambda^o$ increases, the expected injections rise substantially. Although a larger $\lambda^o$ makes annuities cheaper and lowers the annuitization boundary, it also raises the switching level $a_\beta$, expanding the range of annuity income for which capital support remains optimal. Hence, a reduction in annuity price does not necessarily reduce the expected cost of the solvency guarantee.

\section{Conclusion}\label{sec:conclusion}
In this paper, we study the optimal investment, consumption, and annuitization decisions of a pension-plan manager who seeks to balance retirement consumption benefits against the cost of maintaining the pension account's solvency. The pension sponsor provides costly capital injections whenever necessary to prevent negative liquid wealth, while annuity purchases are irreversible. Using the Legendre transform, we derive a semi-explicit solution to the associated free-boundary problem and verify the optimality of the resulting strategies. We identify an endogenous threshold separating injection-active and injection-inactive regions and characterize six policy regimes arising from the interplay between capital injections and consumption constraints. We also obtain a different solution structure in the zero-interest-rate case.

Our results highlight the significant impact of solvency guarantees on optimal pension management. Contrary to the classical benchmark, capital injections in this work delay annuitization and allow positive risky investment even at zero liquid wealth.  The numerical analysis further reveals that cheaper annuities may increase the expected capital contributions required from the pension sponsor. These findings uncover interesting discrepancies in joint decisions of annuitization, investment, and consumption under the funding costs of solvency guarantees in pension-plan management. Future extensions may incorporate stochastic mortality, stochastic interest rates, or surrenderable annuities.

\ \\
\ \\
\textbf{Acknowledgements}: Xiaoqing Liang is supported by the National Natural Science Foundation of China (No. 12371468 and 12271274), and the Natural Science Foundation of Hebei Province (No. A2025202020). Jie Xiong is supported by the National Key R\&D Program of China (No. 2022YFA1006102) and the National Natural Science Foundation of China (No. 12471418).  Xiang Yu is supported by the Hong Kong RGC General Research Fund (GRF) (No. 15304122, 15306523, and 15214125).  \\

\appendix
\section{Special Case with $r=0$}\label{sec:special-cases}
This section presents the result for the special case when the interest rate $r=0$. When the interest rate  $r=0$, the characteristic root $\nu_+$ equals one, and the particular solutions contain logarithmic terms.
We derive these solutions below, distinguishing between the
no-injection endpoint $y_0^{\rm N}(a)$ and the endpoint
$y_0(a)$ when capital injection is available.

\begin{corollary}[Zero-rate specialization]\label{cor:zero-rate-limit}
For $r=0$, define $\theta_0=1/\lambda^o$ and $m=\mu^2/(2\sigma^2)$.
The characteristic roots are $\nu_+=1$ and $\nu_-=-\delta/m$. The coefficients $D_j(a)$ are determined by the
zero-wealth boundary condition \eqref{eq:dual-boundaries} and the
annuitization condition
\eqref{eq:dual-boundaries-annuitization-condition1},
with $\theta=\theta_0$.
The additional condition is \eqref{eq:dual-state-constraint}
when injection is inactive and $y_0(a)=\beta$ when it is active.
At $y=1$, the values and first $y$-derivatives of the two branches match.
The resulting dual solution is
\begin{align}
 \widehat v_+(y,a)
 ={}&D_1(a)y+D_2(a)y^{-\delta/m}
 -\frac{a}{m+\delta}y\log y-\frac{y(\log y)^2}{2\gamma(m+\delta)}
 +\frac{2m+\delta}{\gamma(m+\delta)^2}y\log y,
 \qquad 0<y\leq1, \label{eq:dual-general-low} \\
 \widehat v_0(y,a)
 ={}&D_3(a)y+D_4(a)y^{-\delta/m}
 -\frac{a}{m+\delta}y\log y-\frac1{\gamma\delta},
 \qquad y\geq1.
 \label{eq:dual-general-high}
\end{align}
The boundary ratio $z_{0b}>1$ is the unique solution of
\begin{align}\label{eq:z0b-zero}
 \theta_0(m+\delta)-\log z_{0b}
 +\frac{\delta}{m+\delta}
       \left(z_{0b}^{-(m+\delta)/m}-1\right)=0.
\end{align}
Define
\begin{align}\label{eq:eta}
 \eta=\frac{m+\delta z_{0b}^{-(m+\delta)/m}}{m+\delta}
     =1-\theta_0(m+\delta)+\log z_{0b}\in(0,1).
\end{align}
With $\eta$ defined in \eqref{eq:eta}, the no-injection endpoint is
\begin{align}\label{eq:y0-zero}
 y_0^{\rm N}(a)=
 \begin{cases}
 \left(\dfrac{\eta}{\gamma a}\right)^\eta,
       &0<a\leq\eta/\gamma,\\[0.7em]
 e^{\eta-\gamma a},&a\geq\eta/\gamma.
 \end{cases}
\end{align}
With capital injection, the endpoints are
$y_0(a)=\min\{\beta,y_0^{\rm N}(a)\}$ and
$y_b(a)=y_0(a)/z_{0b}$.
In this case, $a_2=\eta/\gamma$.
The thresholds $a_1$ and $a_\beta$ are given by
$a_1=\kappa z_{0b}^{-1/\kappa}/\gamma$ and \eqref{eq:abeta-positive},
respectively, with $\eta$ replacing $\kappa$.
The six-case classification, the inverse Legendre transform,
and the dual feedback rules remain valid, with
\eqref{eq:dual-general-low}--\eqref{eq:dual-general-high} replacing
the results when the interest rate is positive.
\end{corollary}
\begin{proof}
At $r=0$, the characteristic polynomial factors as
$(\nu-1)(m\nu+\delta)$, yielding the roots $\nu_+$ and $\nu_-$.
The homogeneous differential operator maps $y\log y$ to
$(m+\delta)y$ and $y(\log y)^2$ to
$2(m+\delta)y\log y+2my$. Substituting these identities into the dual
equation \eqref{eq:dual-ode} with $r=0$ yields the particular solutions
in \eqref{eq:dual-general-low}--\eqref{eq:dual-general-high}.
Substituting both branches into the annuitization condition
\eqref{eq:dual-boundaries-annuitization-condition1},
with $\theta=\theta_0$, and solving the resulting equations for
the derivatives of the coefficients yield
\begin{align}
 D_1'(a)&=\theta_0+\frac{\log y_b(a)}{m+\delta}
                  +\frac{m}{(m+\delta)^2},\label{eq:zero-rate-coefficient-derivatives-D1}\qquad
 D_2'(a)=-\frac{m}{(m+\delta)^2}y_b(a)^{1+\delta/m}. 
\end{align}
Matching the values and first derivatives of the expressions
\eqref{eq:dual-general-low} and \eqref{eq:dual-general-high}
at $y=1$ yields coefficient differences that are independent of $a$.
Hence, $D_3'(a)=D_1'(a)$ and $D_4'(a)=D_2'(a)$.
Substituting \eqref{eq:zero-rate-coefficient-derivatives-D1} into
$\widehat v_{ay}$ 
yields \eqref{eq:z0b-zero}. The left-hand side is positive at one,
tends to $-\infty$ as its argument tends to infinity, and has derivative
\[
 -\frac1z-\frac{\delta}{m}z^{-(m+\delta)/m-1}<0,
\]
so the boundary ratio is unique, and \eqref{eq:eta} implies
$0<\eta<1$.
Differentiating \eqref{eq:dual-general-low}--\eqref{eq:dual-general-high}
once with respect to $a$ and twice with respect to $y$, then substituting
\eqref{eq:zero-rate-coefficient-derivatives-D1} and \eqref{eq:eta}, yields
$\widehat v_{ayy}(y_0(a),a)=-\eta/(m y_0(a))$.
Next, differentiate the super-contact condition in
\eqref{eq:positive-rate-contact} along $y=y_0^{\rm N}(a)$
and use the $y$-derivative of \eqref{eq:dual-ode} with $r=0$.
This gives the following equation for the endpoint when injection
is inactive:
\[
 \left(a+\widetilde U'(y_0^{\rm N}(a))\right)
 \frac{(y_0^{\rm N})'(a)}{y_0^{\rm N}(a)}=-\eta.
\]
Integrating this equation on the two consumption branches,
as in \eqref{eq:positive-endpoint-integrated} and
\eqref{eq:mixed-endpoint-solution}, and imposing the bounds
$-e^{-\gamma a}/(\gamma\delta)\leq v^{\rm N}(0,a)\leq0$
yields $a_2=\eta/\gamma$ and \eqref{eq:y0-zero}.

We now determine the coefficients $D_j(a)$ in the region
where capital injection is inactive, $a\geq a_\beta$.
If $y_0^{\rm N}(a)\leq1$, equivalently $a\geq a_2$, the endpoint
lies on the positive-consumption branch. Substituting
\eqref{eq:dual-general-low} into \eqref{eq:positive-rate-contact} gives
\begin{align}
 D_1(a)
 ={}&\frac{a}{m+\delta}
      \left(\log y_0^{\rm N}(a)+1+\frac{m}{m+\delta}\right)
 +\frac{\bigl(\log y_0^{\rm N}(a)\bigr)^2}{2\gamma(m+\delta)}
   -\frac{3m^2+3m\delta+\delta^2}{\gamma(m+\delta)^3},
   \label{eq:zero-rate-inactive-D1}\\
 D_2(a)
 ={}&\frac{m^2\bigl(y_0^{\rm N}(a)\bigr)^{1+\delta/m}}
           {\delta(m+\delta)^2}
   \left(a+\frac{\log y_0^{\rm N}(a)}{\gamma}
            -\frac{m}{\gamma(m+\delta)}\right).
   \label{eq:zero-rate-inactive-D2}
\end{align}
If $y_0^{\rm N}(a)\geq1$, equivalently $0<a\leq a_2$,
the endpoint lies on the zero-consumption branch.
Applying the same two conditions given in \eqref{eq:positive-rate-contact} to
\eqref{eq:dual-general-high} yields
\begin{align}
 D_3(a)
 &=\frac{a}{m+\delta}
   \left(\log y_0^{\rm N}(a)+1+\frac{m}{m+\delta}\right),
   \label{eq:zero-rate-inactive-D3}\\
 D_4(a)
 &=\frac{m^2a}{\delta(m+\delta)^2}
   \bigl(y_0^{\rm N}(a)\bigr)^{1+\delta/m}.
   \label{eq:zero-rate-inactive-D4}
\end{align}
These coefficients determine the entire continuation solution
when $y_0^{\rm N}(a)\geq z_{0b}$, equivalently $0<a\leq a_1$.
In the mixed case $1<y_0^{\rm N}(a)<z_{0b}$,
equivalently $a_1<a<a_2$, the positive-consumption coefficients
are obtained by matching the values and first $y$-derivatives
at $y=1$. Specifically,
\begin{align*}
 (D_1(a)-D_3(a))+(D_2(a)-D_4(a))
 &=-\frac1{\gamma\delta},\\
 (D_1(a)-D_3(a))-\frac{\delta}{m}(D_2(a)-D_4(a))
 &=-\frac{2m+\delta}{\gamma(m+\delta)^2}.
\end{align*}
Solving this system gives
\begin{align*}  
 D_1(a)
 &=D_3(a)-\frac{3m^2+3m\delta+\delta^2}{\gamma(m+\delta)^3},\qquad
 D_2(a)
 =D_4(a)-\frac{m^3}{\gamma\delta(m+\delta)^3}.
\end{align*}
At $y_0^{\rm N}(a)=1$, these expressions agree with
\eqref{eq:zero-rate-inactive-D1}--\eqref{eq:zero-rate-inactive-D2}.
Differentiating the coefficients, using 
\eqref{eq:eta}, \eqref{eq:y0-zero}, and
$y_b^{\rm N}(a)=y_0^{\rm N}(a)/z_{0b}$, recovers
\eqref{eq:zero-rate-coefficient-derivatives-D1}.
Thus, the annuitization condition
\eqref{eq:dual-boundaries-annuitization-condition1} is also satisfied.

For $a<a_\beta$, set $y_0(a)=\beta$ and integrate
\eqref{eq:zero-rate-coefficient-derivatives-D1}  with respect to $a$,
using the values $D_j(a_\beta)$ obtained from the corresponding
solution in the injection-inactive region as initial conditions. The remaining coefficients in the
mixed case follow from the same matching relations at $y=1$. \footnote{These expressions can also be obtained by taking the limit
$r\downarrow0$. We omit the details of this limiting argument for brevity.}

\end{proof}
\begin{remark}[Discounted Shortfall Representation]\label{remark:shortfall}
Let us introduce $\widetilde X_t:=e^{-rt}X_t$, $\widetilde L_t:=\int_0^t e^{-rs}dL_s^X$ and define the unsupported discounted wealth path
\begin{align*}
\widetilde Y_t:=x
+\int_0^t e^{-rs}\big((\mu-r)\pi_s-c_s+A_s\big)ds
+\int_0^t e^{-rs}\sigma\pi_sdW_s
-\int_0^t e^{-rs}\theta dA_s.
\end{align*}
Thus, $\widetilde Y$ represents the discounted wealth trajectory that would result from the same consumption, investment, and annuitization decisions if no capital support were provided.
Integration by parts gives $\widetilde X_t=\widetilde Y_t+\widetilde L_t$ and  the Skorokhod  reflection map yields
\begin{align}\label{eq:minimal-injection}
\quad
\widetilde L_t
=\left(-\inf_{0\leq s\leq t}\widetilde Y_s\right)^+
=0\vee\sup_{0\leq s\leq t}(-\widetilde Y_s),\qquad t\geq0.
\quad
\end{align}
Hence the discounted cumulative capital injection coincides exactly with the running maximum shortfall of the discounted unsupported wealth process. 
Equation \eqref{eq:minimal-injection} is the exact positive-interest-rate analogue of the zero-rate representation in \citet{BoHuangYu25}.  It also clarifies an important point: one can not simply define a candidate unsupported wealth process via   $Y_t:= \int_0^t r\widetilde Y_sds +\int_0^t e^{rs} d \widetilde Y_s$
and expect the original additive decomposition $X_t=Y_t+L_t$ to remain valid. The difficulty arises because the capital injections, once credited to the wealth account, themselves earn the risk-free return.  In contrast, working directly with the original reflected dynamics \eqref{state-X}, or equivalently with the discounted pair $(\widetilde X,\widetilde L)$, preserves a two-dimensional Markov state.
\end{remark}

\section{Proofs}\label{sec:proof}

This section collects all the proofs of the main results in the paper.

\begin{proof}[Proof of Lemma \ref{lem:boundary-ratio}]
The expression in square brackets in \eqref{eq:z0b-general} equals one
at $z_{0b}=1$.  Its derivative, with $z$ denoting the argument, is
\begin{align}\label{eq:positive-rate-root-monotonicity}
 \frac{(\nu_+-1)(1-\nu_-)}{\nu_+-\nu_-}
 \left(\nu_+z^{\nu_+-2}-\nu_-z^{\nu_--2}\right)>0,\qquad z>1.
\end{align}
It tends to infinity as $z\to\infty$, whereas $0<1-r\theta<1$.
The intermediate value theorem and strict monotonicity give the existence
and uniqueness.  Finally,
$0<(1-r\theta)z_{0b}^{\nu_--1}<1$, together with the substitution in the definition of
$\kappa$, results in the desired bounds.
\end{proof}

\begin{proof}[Proof of Lemma \ref{lem:positive-inactive-signs}]
In the positive-consumption case,
\eqref{eq:inactive-positive-coefficientC1}
gives $C_1^{\rm N}(a)<0<C_2^{\rm N}(a)$. The zero-consumption coefficients in
\eqref{eq:inactive-zero-coefficient-C3} satisfy $C_3^{\rm N}(a)<0<C_4^{\rm N}(a)$.
In the mixed case, \eqref{eq:inactive-mixed-coefficients} and
$\Delta_+<0$ from \eqref{eq:general-r-matching-data+} give
$C_1^{\rm N}(a)<0$.
Using \eqref{eq:inactive-zero-coefficient-C3} and
\eqref{eq:general-r-matching-data+}, we also obtain
$C_2^{\rm N}(a)>0$, since
\[
 \gamma a y_0^{\rm N}(a)^{1-\nu_-}
 =\kappa y_0^{\rm N}(a)^{1-\nu_--1/\kappa}
 >\frac1{1-\nu_-},
\]
by \eqref{eq:mixed-endpoint-solution},
$y_0^{\rm N}(a)>1$, and the bound
$\kappa>1/(1-\nu_-)$.
Consequently, on the positive-consumption branch,
\[
 h_y^{\rm N}(y,a)
 =\nu_+(\nu_+-1)^2C_1^{\rm N}(a)y^{\nu_+-2}
  +\nu_-(\nu_--1)^2C_2^{\rm N}(a)y^{\nu_--2}<0.
\]
The same calculation with $C_3^{\rm N}(a)$ and $C_4^{\rm N}(a)$
applies on the zero-consumption branch, where
the relations \eqref{eq:general-r-matching} give value and slope
matching at $y=1$. Evaluating \eqref{eq:dual-ode} on either
side, and then doing the same for its $y$-derivative, gives
matching of $\widehat v_{yy}^{\rm N}$ and
$\widehat v_{yyy}^{\rm N}$, since $\widetilde U$ and
$\widetilde U'$ are continuous by \eqref{eq:utility-conjugate}.
As a result, $h^{\rm N}$ and $h_y^{\rm N}$ are continuous at $y=1$.
Noting $h^{\rm N}(y_0^{\rm N}(a),a)=0$, we obtain
$h^{\rm N}>0$ for $y<y_0^{\rm N}(a)$ and hence strict dual convexity.
The zero-wealth condition in \eqref{eq:positive-rate-contact}
then leads to $\widehat v_y^{\rm N}(y,a)<0$ in the interior.

Recall $g(y,a)=\widehat v_a(y,a)-\theta y$, the expression \eqref{eq:general-r-g-explicit} gives $g_{yy}<0$. The function $\widehat v_{ay}^{\rm N}(y,a)$ equals $\theta$ at
$y=y_b^{\rm N}(a)$ and zero at $y=y_0^{\rm N}(a)$.
It therefore satisfies
$0\leq\widehat v_{ay}^{\rm N}\leq\theta$
and $0\leq\widehat v_{a}^{\rm N}\leq\theta y$.
Meanwhile, we have $b^{\rm N}(a)=-\widehat v_y^{\rm N}(y_b^{\rm N}(a),a) > 0$. Using
\eqref{eq:dual-boundaries-annuitization-condition1}, \eqref{eq:positive-endpoint-solution}, and \eqref{eq:mixed-endpoint-solution} gives
\begin{align*}
 (b^{\rm N})'(a)&=-\theta-(y_b^{\rm N})'(a)
       \widehat v_{yy}^{\rm N}(y_b^{\rm N}(a),a)>-\theta.
\end{align*}
The zero-wealth controls in \eqref{eq:inactive-c-at-zero}, together
with zero consumption when $a\leq a_2$, give
$\pi^{\rm N}(0,a)=0$ and $c^{\rm N}(0,a)=(a-a_2)^+$.
\end{proof}

\begin{proof}[Proof of Proposition \ref{prop:dual-primal}]
The derivative $\widehat v_y(y,a)+x$ increases from $x-b(a)$ to $x$,
so its unique zero is $y^*(x,a)$.
The Legendre identities \eqref{eq:dual-relations} give
$v_x=y^*>0$, $v_{xx}=-1/\widehat v_{yy}<0$, and
$v_a=\widehat v_a>0$, with the dual derivatives evaluated at
$(y^*(x,a),a)$. Substituting these dual expressions into
\eqref{eq:primal-continuation} and simplifying yields
\eqref{eq:dual-ode}.
At zero wealth, $y^*(0,a)=y_0(a)\leq\beta$.
The annuitization condition
\eqref{eq:dual-boundaries-annuitization-condition1}  gives
\eqref{eq:primal-smooth-fit}, since
\[
 v_a(b(a),a)=\widehat v_a(y_b(a),a)=\theta y_b(a),\qquad
 v_{ax}(b(a),a)
 =-\frac{\widehat v_{ay}(y_b(a),a)}{\widehat v_{yy}(y_b(a),a)}
 =\theta v_{xx}(b(a),a).
\]
Finally, the inequality
$\widehat v_a(y^*(x,a),a)-\theta y^*(x,a)\leq0$,
induced by \eqref{eq:general-r-g-explicit}, implies
$v_a(x,a)-\theta v_x(x,a)\leq0$.
These relations hold in the interior, with boundary values
interpreted through the appropriate one-sided limits.
\end{proof}

\begin{proof}[Proof of Lemma \ref{lem:annuity-structure}]
Consider first the injection-active region $0<a<a_\beta$. By \eqref{eq:positive-rate-capped-endpoints}, $y_b(a)={\beta}/{z_{0b}}$, and hence $y_b'(a)=0$. Using the boundary representation
$b(a)=-\hat v_y(y_b(a),a)$ from \eqref{eq:positive-boundary-recovery-ba}, together with the boundary condition $\hat v_{ay}(y_b(a),a)=\theta$, we obtain
\begin{align*}
G'(a)=b'(a)+\theta=-\frac{d(\hat{v}_{y}(y_b(a),a))}{d a}+\theta=-\hat{v}_{yy}(y_b(a),a)y_b'(a)-\hat{v}_{ay}(y_b(a),a)+\theta=0.
\end{align*}
For $a>a_{\beta}$, Lemma~\ref{lem:positive-inactive-signs} gives $b'(a)>-\theta$, which yields $ G'(a)>0$. Finally, in view that $b(a)\geq0$ and $\theta>0$,
$G(a)=b(a)+\theta a \geq \theta a$, and hence $G(a)\to\infty$ as $a\to\infty$.   This completes the proof. 
\end{proof}

\begin{proof}[Proof of Theorem \ref{thm:positive-rate-candidate}]
Lemma~\ref{lem:positive-inactive-signs} and
\eqref{eq:positive-rate-curvature-comparison} assert the convexity of
the dual candidate. Proposition~\ref{prop:dual-primal} then gives
monotonicity in both state variables, strict concavity in wealth,
the HJB equality in the continuation region, and the boundary
conditions. It remains to verify the HJB-VI \eqref{eq:HJB-v} in
the annuitization region, smooth fit and the boundedness.

For $x>b(a)$, the immediate annuity purchase moves the state to
$(b(\bar a),\bar a)$, where $\bar a>a_\beta$ and capital injection
is inactive. Hence $b'(\bar a)+\theta>0$. Differentiating
$x+\theta a=b(\bar a)+\theta\bar a$ and
\eqref{eq:annuitization-value-extension} gives
\begin{align*}
 v_x(x,a)=y_b(\bar a),\qquad
 v_a(x,a)=\theta y_b(\bar a),\qquad
 v_{xx}(x,a)=\frac{y_b'(\bar a)}{b'(\bar a)+\theta}
           =-\frac1{\widehat v_{yy}(y_b(\bar a),\bar a)}.
\end{align*}
The last equality follows by differentiating the boundary relation
\eqref{eq:positive-boundary-recovery-ba} at $\bar a$ and applying
the annuitization condition $\widehat v_{ay}(y_b(\bar{a}),\bar{a})=\theta.$
These identities also establish monotonicity in both state
variables and strict concavity in wealth in the annuitization region.
Moreover, $v$, $v_x$, and $v_{xx}$ at $(x,a)$ coincide with their
values at $(b(\bar a),\bar a)$. The optimized investment and
consumption terms in the Hamiltonian are therefore identical at
the two points. The difference between the Hamiltonians is
\begin{align}\label{eq:general-r-action-drift}
 \bigl[(rx+a)-(r\bar x+\bar a)\bigr]y_b(\bar a)
 =-(1-r\theta)(\bar a-a)y_b(\bar a)\leq0.
\end{align}
At the boundary point, the Hamiltonian is zero by
\eqref{eq:primal-continuation}. Equation~\eqref{eq:general-r-action-drift}
therefore proves \eqref{eq:HJB-v} in the annuitization region, where $v_a-\theta v_x=0$ also holds.

The preceding identities establish smooth fit along the
annuitization boundary where injection is inactive. Where
injection is active, both $v(b(a),a)$ and $y_b(a)=\beta/z_{0b}$
are constant. The value function and its first derivatives
therefore match across this boundary segment as well.
Second derivatives with respect to wealth need not agree across
the active segment\footnote{The annuitization super-contact condition
is $v_{ax}=\theta v_{xx}$, interpreted from the continuation region.
Along $x=b(a)$, it implies
$y_b'(a)=(b'(a)+\theta)v_{xx}$.
For $a<a_\beta$, we have $y_b'(a)=0$ and $b'(a)=-\theta$,
so this identity reduces to $0=0$ and does not require the two
one-sided limits of $v_{xx}$ to coincide.}
and may diverge as the zero-wealth boundary
is approached in the inactive region.

To prove the boundedness, substitute
\eqref{eq:inactive-positive-coefficientC1}
and $\log y_0^{\rm N}(a)=\kappa-\gamma a$ into
\eqref{eq:dual-general-r-low}. For $a\geq a_2$, this yields
\begin{align}\label{eq:positive-rate-tail-bound}
 \widehat v^{\rm N}(y,a)
 =y_0^{\rm N}(a)\,
   \widehat v^{\rm N}\left(\frac{y}{y_0^{\rm N}(a)},a_2\right),
 \qquad
 \frac1{z_{0b}}\leq\frac{y}{y_0^{\rm N}(a)}\leq1.
\end{align}
Furthermore, it follows from formula~\eqref{eq:inactive-boundary-bN} that 
$b^{\rm N}(a)=b^{\rm N}(a_2)$, for any $a\geq a_2$.
Thus, the annuitization boundary is constant for $a\geq a_2$.
For $0\leq x\leq b^{\rm N}(a_2)$, applying the inverse Legendre transform
\eqref{eq:inverse-transform} and changing variables from $y$ to
$z=y/s$, where $s:=y_0^{\rm N}(a)=e^{\kappa-\gamma a}=e^{-\gamma(a-a_2)}$, we obtain
\begin{align*}
 v^{\rm N}(x,a)
 &=\inf_{s/z_{0b}\leq y\leq s}
   \left\{\widehat v^{\rm N}(y,a)+xy\right\}
 =s\inf_{1/z_{0b}\leq z\leq1}
   \left\{\widehat v^{\rm N}(z,a_2)+xz\right\}
 =s\,v^{\rm N}(x,a_2).
\end{align*}
For $x>b^{\rm N}(a_2)$, the immediate annuity purchase has size
$\Delta a=(x-b^{\rm N}(a_2))/\theta$. By the same value-matching construction
as in \eqref{eq:annuitization-value-extension}, we have
\begin{align}
 v^{\rm N}(x,a)
 =v^{\rm N}(b^{\rm N}(a_2),a+\Delta a)
 &=e^{-\gamma(a + \Delta a -a_2)}v^{\rm N}(b^{\rm N}(a_2),a_2)\\*
 &=e^{-\gamma(a-a_2)}v^{\rm N}(b^{\rm N}(a_2),a_2+\Delta a)\\*
 &=e^{-\gamma(a-a_2)}v^{\rm N}(x,a_2).\label{eq:positive-rate-primal-scaling}
\end{align}
Here, the second and third equalities follow by applying the scaling at
$x=b^{\rm N}(a_2)$ to the income levels $a+\Delta a$ and $a_2+\Delta a$, respectively.

For $a\geq\max\{a_2,a_\beta\}$, the candidate value function
with capital injection coincides with $v^{\rm N}(x,a)$.
Thus, \eqref{eq:positive-rate-primal-scaling} shows that the
candidate converges exponentially to zero as $a\to\infty$
for each fixed $x$.
On each bounded income interval, its dual endpoints remain
positive and finite, and its coefficient functions extend
continuously to $a=0$. The explicit dual branches and the inverse
transform therefore give a bounded candidate on the corresponding
continuation region. Together with the exponential scaling and
the constant boundary for large $a$, this proves boundedness
throughout the continuation region. Value matching then gives
boundedness in the annuitization region. Finally, monotonicity
in $a$ and the limit $v(x,a)\to0$ as $a\to\infty$ imply $v\leq0$.

Applying the feedback formulas \eqref{eq:dual-feedback-detailed}
to \eqref{eq:positive-rate-tail-bound} also gives, for
$a\geq a_2$ and $0\leq x\leq b(a_2)$,
\begin{align}\label{eq:positive-tail-policies}
 \pi^{\rm N}(x,a)=\pi^{\rm N}(x,a_2),\qquad
 c^{\rm N}(x,a)=a-a_2+c^{\rm N}(x,a_2).
\end{align}
The endpoint values are understood by continuity. Together with
continuity on bounded income intervals, these identities imply
that, throughout the continuation region, risky investment is
bounded and consumption grows at most linearly in $a$. The annuitization boundary
is also bounded as it is continuous on bounded income intervals and constant for sufficiently large $a$.
\end{proof}

\begin{proof}[Proof of Theorem \ref{thm:verification}]
We first show $v(x,a)\geq V(x,a)$ for any $(x,a)\in\mathbb{R}_+^2$. Fix an arbitrary admissible control $(\pi,c,A)$. 
 Define
\begin{align*}
 \mathcal H^{\pi,c}v
 &:=(rx+(\mu-r)\pi-c+a)v_x
   +\frac12\sigma^2\pi^2v_{xx}
   -\delta v+U(c).
\end{align*}
For $n\in\mathbb{N}$, define the stopping time
\begin{align*}
\tau_n:=n\wedge\inf\left\{t\geq0:
 X_t+A_t+L_t+\int_0^t(c_s+\pi_s^2)\,ds\geq n\right\},
\end{align*}
and, for a fixed $T>0$, set $T_n:=T\wedge\tau_n$. Applying the  It\^o formula to $e^{-\delta t}v(X_t,A_t)$ on $[0,T_n]$ gives
\begin{align}\label{eq:ito-verification}
 e^{-\delta \tau_n}v(X_{T_n},A_{T_n})
 =&v(x,a)
 +\int_0^{T_n}e^{-\delta t}
   \big(\mathcal H^{\pi_t,c_t}v-U(c_t)\big)dt \notag\\
 &+\int_0^{T_n}e^{-\delta t}v_x(0,A_t)dL_t^c
 +\int_0^{T_n}e^{-\delta t}
   (v_a-\theta v_x)(X_t,A_t)dA_t^c \notag\\
 &+\sum_{0<t\leq T_n}e^{-\delta t}
   \big[v(X_t,A_t)-v(X_{t-},A_{t-})\big]
 +M_{T_n},
\end{align}
where $A^c$ and $L^c$ denote the continuous parts of $A$ and $L$,
respectively, and
\begin{align*}
 M_t:=\int_0^t e^{-\delta s}\sigma\pi_sv_x(X_{s-},A_{s-})dW_s.
\end{align*}

Since $v$ is non-decreasing and concave in $x$ by Theorem \ref{thm:positive-rate-candidate} and Corollary \ref{cor:zero-rate-limit},  $0\leq v_x(x,a)\leq v_x(0,a)\leq\beta$, for all $(x,a)\in\R_+^2$. Consequently, admissibility of $\pi$ implies that, for every finite $T$,
\begin{align*}
\Ex\left[\int_0^T e^{-2\delta s}\sigma^2\pi_s^2 v_x^2(X_{s-},A_{s-})ds\right]
\leq\sigma^2\beta^2\Ex\left[\int_0^T\pi_s^2ds\right]<\infty,
\end{align*}
Thus $M$ is a square-integrable martingale on every finite horizon, and in particular $\Ex[M_{T_n}]=0$. 

The HJB-VI \eqref{eq:HJB-v} implies, for every admissible $(\pi,c, A)$,
\begin{align}\label{eq:verification-HJB-inequalities}
 \mathcal H^{\pi,c}v\leq0,
 \qquad
 v_a-\theta v_x\leq0.
\end{align}
Because the capital injection process $L$ increases only on $\{X_t=0\}$, the lower-boundary condition $v_x(0,a)\leq\beta$ gives
\begin{align}\label{eq:continuous-injection-bound}
 \int_0^{T_n}e^{-\delta t}
   \big(v_x(0,A_t)-\beta\big)dL_t^c\leq0.
\end{align}
Likewise, the annuitization inequality in
\eqref{eq:verification-HJB-inequalities} yields
\begin{align}\label{eq:continuous-annuity-bound}
\int_0^{T_n}e^{-\delta t}
   (v_a-\theta v_x)(X_t,A_t)dA_t^c \leq 0.
\end{align}

It remains to control the jump terms. This also covers the
possibility that $A$ and $L$ jump simultaneously. Fix a jump time $t$ and write   $x=X_{t-}$, $a=A_{t-}$, $z=\Delta A_t$, and $\ell=\Delta L_t$.  Since
$X_t=x-\theta z+\ell\geq0$, the whole segment
$u\mapsto(x+\ell-\theta u,a+u)$, $0\leq u\leq z$, remains in $\mathbb R_+^2$.  We may therefore decompose the jump increment as
\begin{align}\label{eq:joint-jump-bound}
 &v(X_t,A_t)-v(X_{t-},A_{t-})-\beta\Delta L_t \notag\\
 &\quad=\big[v(x+\ell-\theta z,a+z)
              -v(x+\ell,a)\big]
       +\big[v(x+\ell,a)-v(x,a)-\beta\ell\big]\notag\\
 &\quad=\int_0^z(v_a-\theta v_x)
          (x+\ell-\theta u,a+u)du
       +\int_0^\ell(v_x(x+s,a)-\beta)ds
 \leq0.
\end{align}
The last inequality follows from
\eqref{eq:verification-HJB-inequalities} and
$0\leq v_x\leq\beta$.

By \eqref{eq:verification-HJB-inequalities}-\eqref{eq:joint-jump-bound}, all finite-variation terms on the right-hand side of \eqref{eq:ito-verification} after the term $v(x,a)$ are non-positive. Taking expectations and using $\Ex[M_{T_n}]=0$ therefore yields
\begin{align}\label{eq:finite-horizon-verification}
 &\Ex\left[\int_0^{T_n}e^{-\delta t}U(c_t)dt
      -\beta\int_0^{T_n}e^{-\delta t}dL_t\right]
\leq v(x,a)-\Ex\left[e^{-\delta T_n}
             v(X_{T_n},A_{T_n})\right].
\end{align}
For fixed $T$, admissibility implies $T_n\uparrow T$ almost surely. By Theorem \ref{thm:positive-rate-candidate} and Corollary \ref{cor:zero-rate-limit}, $v$ is bounded, thus, the terminal term converges by dominated convergence.
Then, the  dominated convergence for the utility term and monotone convergence for the injection
term give
\begin{align}\label{eq:finite-T-verification}
 &\Ex\left[\int_0^T e^{-\delta t}U(c_t)dt
 -\beta\int_0^T e^{-\delta t}dL_t\right]\leq v(x,a)-\Ex\left[e^{-\delta T}v(X_T,A_T)\right].
\end{align}
Finally, let $T\to\infty$.  The boundedness of $v$ implies that the terminal expectation in \eqref{eq:finite-T-verification} converges to zero.  Using the dominated convergence theorem  and monotonic convergence theorem once more gives
\begin{align*}
 \Ex\left[\int_0^\infty e^{-\delta t}U(c_t)dt
 -\beta\int_0^\infty e^{-\delta t}dL_t\right]\leq v(x,a).
\end{align*}
Since $(\pi,c,A)$ was arbitrary, taking the supremum over all admissible controls yields $V(x,a)\leq v(x,a)$.

We now show that $v(x,a)= V(x,a)$ for any $(x,a)\in\mathbb{R}_+^2$ and the candidate control $(\pi^*,c^*,A^*)$ is optimal. By the proof of Theorem~\ref{thm:positive-rate-candidate}, we know that $\pi^*(x,a)$ is bounded in $(x,a)$ and that $c^*(x,a)$ grows at most linearly in $a$. Hence, for every $T>0$,
\begin{align*}
\mathbb{E}\left[\int_0^T\left((\pi_t^*)^2+c_t^*\right)dt\right]<\infty.
\end{align*}
which shows that $(\pi^*,c^*,A^*)\in\mathbb{U}$. By \eqref{eq:dual-feedback-detailed}, the regular feedback controls $(\pi^*,c^*)$ attain the continuation Hamiltonian equality, so $
\mathcal H^{\pi_t^*,c_t^*}v (X_t^*,A_t^*)=0$.
  By $(\beta-v_x(0,A_t^*)\big)dL_t^{X^*}=0$, the inequality for capital injection \eqref{eq:continuous-injection-bound} becomes an equality, and thus the capital-injection term in \eqref{eq:ito-verification} vanish. The support condition for $A^{*,c}$ implies that it acts only on the annuitization boundary, where
$v_a-\theta v_x=0$. Hence the continuous-annuitization term also vanishes. 

At a jump time $t$ of $A^*$, let $z=\Delta A_t^*>0$. By the optimal projection rule $X_{t-}^*-\theta z=b(A_t^*)>0$. Hence the annuity purchase alone leaves the wealth process strictly non-negative. Since $L^{X^*}$ is the minimal lower regulator, no simultaneous capital injection is required, and therefore $\Delta L_t^{X^*}=0$, $X_t^*=X_{t-}^*-\theta\Delta A_t^*=b(A_t^*)$. By \eqref{eq:joint-jump-bound} together with $v_a-\theta v_x=0$,  $v(X_t^*,A_t^*) -v(X_{t-}^*,A_{t-}^*) -\beta\Delta L_t^{X^*}=0$, so the jump inequality \eqref{eq:joint-jump-bound} is attained with
equality under the proposed policy.

Consequently, every inequality above becomes an equality under
the candidate control. Passing first through the localization sequence and then letting $T\to\infty$ yields
\begin{align*}
v(x,a)=\Ex\left[\int_0^\infty e^{-\delta t}U(c_t^*)dt
-\beta\int_0^\infty e^{-\delta t}\,dL_t^{X^*}\right]\leq V(x,a).
\end{align*}
Combining this with the established inequality
$V(x,a)\leq v(x,a)$, we get that $v(x,a)=V(x,a)$
and $(\pi^*,c^*,A^*)$ is optimal.
\end{proof}


\begin{thebibliography}{}
{\small
\bibitem[Andr\'easson and Shevchenko(2024)]{AndreassonShevchenko24} Andr\'easson, J. G., Shevchenko, P. V. (2024): Optimal annuitisation, housing and reverse mortgage in retirement in the presence of a means-tested public pension. {\it European Actuarial Journal} {14}: 871--904. 

\bibitem[Blake et al.(2003)]{BlakeCairnsDowd03}
Blake, D., Cairns, A. J. G., Dowd, K. (2003):
Pensionmetrics 2: Stochastic pension plan design during the distribution phase.
{\it Insurance: Mathematics and Economics} {33}(1): 29--47.

\bibitem[Bo et al.(2025a)]{BoHuangYu24} Bo, L., Huang, Y., Yu, X. (2025): Stochastic control problems with state reflections arising from relaxed benchmark tracking. {\it Mathematics of Operations Research} 50(4): 2526-2551.

\bibitem[Bo et al.(2025b)]{BoHuangYu25SICON} Bo, L., Huang, Y., Yu, X. (2025):  On optimal tracking portfolio in incomplete markets: The reinforcement learning approach. {\it SIAM Journal on Control and Optimization} 63(1): 321-348.

\bibitem[Bo et al.(2026a)]{BoHuangYu25} Bo, L., Huang, Y., Yu, X. (2026): An extended Merton problem with relaxed benchmark tracking. {\it Mathematical Finance} 36(2): 422-448.

\bibitem[Bo et al.(2026b)]{BHYY26} Bo, L., Huang, Y., Yan, K., Yu, X. (2026): Optimal consumption under relaxed benchmark tracking and consumption drawdown constraint. {\it SIAM Journal on Financial Mathematics} 17(1):78-117. 

\bibitem[Bo et al.(2021)]{BoLiaoYu21}  Bo, L., Liao, H., Yu, X. (2021): Optimal tracking portfolio with a ratcheting capital benchmark. {\it SIAM Journal on Control and Optimization} {59}(3): 2346-2380.


\bibitem[Browne(1999)]{Browne99} Browne, S. (1999): Beating a moving target: Optimal portfolio strategies for outperforming a stochastic benchmark. {\it Finance and Stochastics} {3}(3): 275-294.

\bibitem[Costantini(1992)]{Costantini1992}
Costantini, C. (1992): The Skorohod oblique reflection problem in domains with corners and application to stochastic differential equations,
{\it Probability Theory and Related Fields}
{91}(1): 43--70.

\bibitem[Cox and Huang(1989)]{CoxHuang89} Cox, J. C., Huang, C.-F. (1989): Optimal consumption and portfolio policies when asset prices follow a diffusion process. {\it Journal of Economic Theory} {49}(1): 33-83.

\bibitem[Dadashi(2020)]{Dadashi20}
Dadashi, H. (2020):
Optimal investment--consumption problem: Post-retirement with minimum guarantee.
{\it Insurance: Mathematics and Economics} {94}: 160--181.

\bibitem[Davidoff et al.(2005)]{DavidoffBrownDiamond05} Davidoff, T., Brown, J. R., Diamond, P. A. (2005): Annuities and individual welfare. {\it American Economic Review} {95}(5): 1573-1590.
%
%


\bibitem[{Deng et al.(2022)}]{DLPY22} Deng, S., Li, X., Pham, H., Yu, X. (2022):  Optimal consumption with reference to past spending maximum. {\it Finance and Stochastics} {26}: 217-266.



\bibitem[Di Giacinto et al.(2014)]{DiGiacintoFedericoGozziVigna14}
Di Giacinto, M., Federico, S., Gozzi, F., Vigna, E. (2014):
Income drawdown option with minimum guarantee.
{\it European Journal of Operational Research} {234}(3): 610--624.

\bibitem[Dillschneider et al.(2024)]{DillschneiderEtAl24} Dillschneider, Y., Maurer, R., Schober, P. (2024): Optimal annuitization under stochastic interest rates. {\it ASTIN Bulletin} {54}(3): 626--651. 

\bibitem[Eisenberg and Schmidli(2009)]{EisenbergSchmidli09} Eisenberg, J., Schmidli, H. (2009): Optimal control of capital injections by reinsurance in a diffusion approximation. {\it Bl\"atter DGVFM} {30}: 1-13.

\bibitem[Eisenberg and Schmidli(2011)]{EisenbergSchmidli11} Eisenberg, J., Schmidli, H. (2011): Optimal control of capital injections by reinsurance with a constant rate of interest. {\it Journal of Applied Probability} {48}(3): 733--748. 

\bibitem[Eisenberg and Kr\"uhner(2018)]{EisenbergKruhner18} Eisenberg, J., Kr\"uhner, P. (2018): The impact of negative interest rates on optimal capital injections. {\it Insurance: Mathematics and Economics} {82}: 1--10. 

\bibitem[Ferrari and Schuhmann(2019)]{FerrariSchuhmann19} Ferrari, G., Schuhmann, P. (2019): An optimal dividend problem with capital injections over a finite horizon. {\it SIAM Journal on Control and Optimization} {57}(4): 2686-2719.

\bibitem[Gaivoronski et al.(2005)]{GaivoronskiEtAl05} Gaivoronski, A., Krylov, S., van der Wijst, N. (2005): Optimal portfolio selection and dynamic benchmark tracking. {\it European Journal of Operational Research} {163}(1): 115-131.

\bibitem[Gerrard et al.(2004)]{GerrardHabermanVigna04}
Gerrard, R., Haberman, S., Vigna, E. (2004):
Optimal investment choices post-retirement in a defined contribution pension scheme.
{\it Insurance: Mathematics and Economics} {35}(2): 321--342.

\bibitem[Gerrard et al.(2012)]{GerrardHojgaardVigna12}
Gerrard, R., H{\o}jgaard, B., Vigna, E. (2012):
Choosing the optimal annuitization time post-retirement.
{\it Quantitative Finance} {12}(7): 1143--1159.

\bibitem[Guo et al.(2026)]{GLSX2026} Guo, J., Liang, X., Shen, Y.,  Xiong, J. (2026): Optimal annuitization and asset allocation with fixed transaction costs. {\it Insurance: Mathematics and Economics} 127:  103223.

\bibitem[Habib et al.(2020)]{HabibEtAl20} Habib, F., Huang, H., Mauskopf, A., Nikolic, B., Salisbury, T. S. (2020): Optimal allocation to deferred income annuities. {\it Insurance: Mathematics and Economics} {90}: 94--104. 

\bibitem[Horneff et al.(2008)]{HorneffMaurerStamos08}
Horneff, W. J., Maurer, R. H., Stamos, M. Z. (2008):
Optimal gradual annuitization: Quantifying the costs of switching to annuities.
{\it Journal of Risk and Insurance} {75}(4): 1019--1038.

\bibitem[Karatzas et al.(1986)]{KaratzasEtAl86} Karatzas, I., Lehoczky, J. P., Sethi, S. P., Shreve, S. E. (1986): Explicit solution of a general consumption/investment problem. {\it Mathematics of Operations Research} {11}(2): 261-294.

\bibitem[Liang and Young(2018)]{LiangYoung18} Liang, X., Young, V. R. (2018): Annuitization and asset allocation under exponential utility. {\it Insurance: Mathematics and Economics} {79}: 167-183.

\bibitem[Lions and Sznitman(1984)]{LionsSznitman84} Lions, P.-L., Sznitman, A.-S. (1984): Stochastic differential equations with reflecting boundary conditions. {\it Communications on Pure and Applied Mathematics} {37}(4): 511--537. 

\bibitem[Lokka and Zervos(2008)]{LokkaZervos08} Lokka, A., Zervos, M. (2008): Optimal dividend and issuance of equity policies in the presence of proportional costs. {\it Insurance: Mathematics and Economics} {42}(3): 954-961.

\bibitem[Merton(1969)]{Merton69} Merton, R. C. (1969): Lifetime portfolio selection under uncertainty: The continuous-time case. {\it Review of Economics and Statistics} {51}(3): 247-257.

\bibitem[Merton(1971)]{Merton71} Merton, R. C. (1971): Optimum consumption and portfolio rules in a continuous-time model. {\it Journal of Economic Theory} {3}(4): 373-413.

\bibitem[Milevsky et al.(2006)]{MilevskyMooreYoung06} Milevsky, M. A., Moore, K. S., Young, V. R. (2006): Asset allocation and annuity-purchase strategies to minimize the probability of financial ruin. {\it Mathematical Finance} {16}(4): 647-671.

\bibitem[Milevsky and Young(2007)]{MilevskyYoung07} Milevsky, M. A., Young, V. R. (2007): Annuitization and asset allocation. {\it Journal of Economic Dynamics and Control} {31}: 3138-3177.

\bibitem[Rauh(2006)]{Rauh06}
Rauh, J. D. (2006):
Investment and financing constraints: Evidence from the funding of corporate pension plans.
{\it Journal of Finance} {61}(1): 33--71.

\bibitem[Saisho(1987)]{Saisho87} Saisho, Y. (1987): Stochastic differential equations for multi-dimensional domain with reflecting boundary. {\it Probability Theory and Related Fields} {74}(3): 455--477. 

\bibitem[Shreve et al.(1984)]{ShreveLehoczkyGaver84}
Shreve, S. E., Lehoczky, J. P., Gaver, D. P. (1984):
Optimal consumption for general diffusions with absorbing and reflecting barriers.
{\it SIAM Journal on Control and Optimization} {22}(1): 55--75.

\bibitem[Tanaka(1979)]{Tanaka1979}
Tanaka, H. (1979): Stochastic differential equations with reflecting boundary condition in convex regions.
{\it Hiroshima Mathematical Journal} {9}: 163--177.

\bibitem[Wang and Young(2012)]{WangYoung12} Wang, T., Young, V. R. (2012): Maximizing the utility of consumption with commutable annuities. {\it Insurance: Mathematics and Economics} {51}(2): 352-369.

\bibitem[Yaari(1965)]{Yaari65} Yaari, M. E. (1965): Uncertain lifetime, life insurance, and the theory of the consumer. {\it Review of Economic Studies} {32}(2): 137-150.

\bibitem[Yin and Yuen(2015)]{YinYuen15} Yin, C., Yuen, K. C. (2015): Optimal dividend problems for a jump-diffusion model with capital injections and proportional transaction costs. {\it Journal of Industrial and Management Optimization} {11}(4): 1247--1262. 
}
\end{thebibliography}
\end{document}